\RequirePackage{fix-cm}

\documentclass[preprint]{imsart}

\RequirePackage[OT1]{fontenc}
\RequirePackage{amsthm,amsmath}
\RequirePackage[numbers]{natbib}
\RequirePackage[colorlinks,citecolor=blue,urlcolor=blue]{hyperref}

\arxiv{arXiv:0000.0000}

\usepackage{float}
\usepackage{breakcites}
\usepackage{epsfig}
\usepackage{amssymb}
\usepackage{amsmath}
\usepackage{hyperref}
\usepackage{verbatim}
\usepackage{subcaption}
\usepackage{tikz}
\usepackage{color}

\usepackage{bbm}

\usetikzlibrary{arrows}
\usepackage[bottom=0.95in,top=0.95 in]{geometry}
\tikzset{cross/.style={cross out, draw=black, minimum size=2*(#1-\pgflinewidth), inner sep=0pt, outer sep=0pt},
cross/.default={5pt}}
\usetikzlibrary{shapes.misc}
\newcommand{\be}{\begin{eqnarray}}
\newcommand{\ee}{\end{eqnarray}}
\newcommand{\bea}{\begin{eqnarray*}}
\newcommand{\eea}{\end{eqnarray*}}

\usepackage{nicefrac}
\numberwithin{equation}{section}

\usepackage{pst-node,pst-plot}
\pstVerb{realtime srand}
\psset{plotpoints=50}

\startlocaldefs
\numberwithin{equation}{section}
\theoremstyle{plain}
\newtheorem{theorem}{Theorem}[section]
\newtheorem{lemma}{Lemma}[section]
\newtheorem{corollary}{Corollary}[section]
\newtheorem{proposition}{Proposition}[section]
\endlocaldefs

\begin{document}

\begin{frontmatter}
\title{First passage time  for Slepian process with linear barrier}
\runtitle{First passage time  for Slepian process with linear barrier}

\begin{aug}
\author{\fnms{Jack} \snm{Noonan}\ead[label=e1]{Noonanj1@cf.ac.uk}},
\author{\fnms{Anatoly} \snm{Zhigljavsky}\ead[label=e2]{ZhigljavskyAA@cardiff.ac.uk}}

\runauthor{J. Noonan and A. Zhigljavsky}

\affiliation{Cardiff University}

\address{School of Mathematics,\\
Cardiff University, Cardiff,\\
CF24 4AG, UK\\
\printead{e1}\\
\phantom{E-mail:\ }\printead*{e2}}
\end{aug}



\begin{abstract}
In this paper we extend results of L.A. Shepp by finding  explicit formulas for the first passage probability $F_{a,b}(T\, |\, x)={\rm Pr}(S(t)<a+bt \text{  for all  } t\in[0,T]\,\, | \,\,S(0)=x)$, for all $T>0$, where $S(t)$ is a Gaussian process with mean 0 and covariance $\mathbb{E} S(t)S(t')=\max\{0,1-|t-t'|\}\,.$  We then extend the results to the case of piecewise-linear barriers and outline applications to change-point detection problems. Previously, explicit formulas for $F_{a,b}(T\, |\, x)$ were  known only  for the cases $b=0$ (constant barrier) or $T\leq 1$ (short interval).
 \end{abstract}

\end{frontmatter}

%
%
%
%
%
%
%
%

\section{Introduction}

Let $T>0$ be a fixed real number and let $S(t)$,  $ t\in [0,T]$, be a Gaussian process with mean 0 and covariance
$$
\mathbb{E} S(t)S(t')=\max\{0,1-|t-t'|\}\, .
$$
This process is often called Slepian process and can be expressed in terms of the standard Brownian motion $W(t)$ by
\be\label{BM}
S(t)=W(t)-W(t+1), \,\, t\ge 0.
\ee
Let $a$ and $b$ be fixed real numbers and $x<a$. We are interested in an explicit formula for the first passage probability
\be\label{problem}
F_{a,b}(T\, |\, x) := {\rm Pr}(S(t)<a+bt \text{  for all   } t\in[0,T]\,\, | \,\,S(0)=x);
\ee
note   $F_{a,b}(T\, |\, x)=0$ for $x \geq a$.

The case of a constant barrier, when $b=0$, has attracted significant attention in literature. In his seminal paper \cite{slepian1961first}, D.Slepian has shown how to derive an explicit expression for $F_{a,0}(T\, |\, x)$ in the case $T\leq 1$; see also \cite{Mehr}. The case $T>1$ is much more complicated than the case $T \leq 1$. Explicit formulas for $F_{a,0}(T\, |\, x)$ with general $T$  were derived by L. A. Shepp in \cite{Shepp71}; these formulas  are  special cases of results formulated in Section~\ref{main_result_1} and \ref{main_result_2}.
We believe our paper can be considered as a natural extension of the methodology developed in \cite{slepian1961first} and \cite{Shepp71}; hence the title of this paper.

In the case $T\leq 1$, Slepian's method for deriving formulas for $F_{a,0}(T\, |\, x)$ can be easily extended to the case of a general linear barrier; see
 Section~\ref{sec:T1} for the discussion and formulas for $F_{a,b}(T\, |\, x)$ with $T\le1$. For general $T>0$, including the case $T>1$,
  explicit formulas for $F_{a,b}(T\, |\, x)$ were unknown. Derivation of these formulas is the main objective of this paper.

To do this,  we  generalise Shepp's methodology of \cite{Shepp71}.
  The principal distinction between Shepp's methodology and our results is the use of
    an alternative way of computing coincidence probabilities.  Shepp's proofs  heavily rely on the so-called Karlin-McGregor identity, see \cite{karlin1959coincidence}, but we use a different result formulated and discussed in Section~\ref{sec:aux}.

   The structure of the paper is as follows. In Section~\ref{sec:Int_T}, we derive an expression for $F_{a,b}(T\, |\, x)$ for integer $T$ and in Section~\ref{sec:Non_int_T}
  we extend the results for non-integral $T$.
           In Sections~\ref{sec:Special_case_barrier} and \ref{Two_slope_changes}, we extend the results to the case of piecewise-linear barriers. In Section~\ref{sec:Change_point}, we outline an application to a change-point detection problem; this application was our main motivation for this research. In Appendix A, we discuss formulas for $F_{a,b}(T\, |\, x)$ with $T\le1$ and provide approximations for the ARL (average run length) in a change-point detection procedure. In Appendix B, we give two technical proofs.

\section{Linear barrier $a+bt$ with integral $T$}\label{sec:Int_T}

In this Section, we  derive an explicit formula for the first passage probability $F_{a,b}(T\, |\, x)$ defined in \eqref{problem} under the assumption  that $T$ is a positive integer, $T=n$. First, we formulate and slightly modify a general result from  \cite[p.40]{katori2012reciprocal}.

\subsection{An important auxiliary result}
\label{sec:aux}

\begin{lemma}\label{Theorem_Kat}
 For any $s>0$ and a positive integer $n$, let  $W_i(t)$, $t\in [0,s]$  be $n+1$ independent Brownian Motion processes with drift parameters $\mu_i\in \mathbb{R}$; $i=0,1,\ldots, n$.
 Suppose $a_0<a_1<\ldots<a_n$ and $c_0<c_1<\ldots<c_n$ and let $dc_0,\ldots,dc_n$ be infinitesimal intervals around $c_0,\ldots,c_n$. Construct the vectors $\boldsymbol{\mu}=(\mu_0,\mu_1,\ldots,\mu_n)'$, $\textbf{a}=(a_0,a_1,\ldots,a_n)'$ and $\textbf{c}=(c_0,c_1,\ldots,c_n)'$. Then
\be\label{Katori}
 \lefteqn{{\rm Pr}\big\{ W_0(t)\!<\!W_1(t)\!<\!\cdots\!<\!W_n(t), \;\!\! \forall t \in [0,s],\!\! \; W_i(s)\in dc_i\; (0\leq i \leq  n) \, \big | W_i(0) =a_i\; (0\leq i \leq  n) \big\} }
 \nonumber\\
 && \quad\quad\quad\quad\quad\quad\quad\quad\quad= \exp\left(-\frac{s}{2}|\boldsymbol{\mu}|^2+\boldsymbol{\mu}\cdot(\textbf{c}-\textbf{a})  \right)\det\left[\varphi_{s}({a}_i-{c}_j) \right]_{i,j=0}^n dc_0 dc_1\ldots dc_n,\;\;\;\;\;
\ee
where $|\cdot  |$ denotes the Euclidean norm, $\cdot$ denotes the scalar  product and $\varphi_{s}(a-c)dc={\rm Pr}(W(s)\in dc \,|\, W(0)=a)$ is the transition probability for the standard Brownian Motion with no drift,
\be\label{BM_transition}
\varphi_s(z) := \frac{1}{\sqrt{2\pi s} }e^{-z^2/(2s)}\, .
\ee
\end{lemma}
Lemma~\ref{Theorem_Kat} is an extension of the celebrated result of Karlin and McGregor on coincidence probabilities (see \cite{karlin1959coincidence}) when applied specifically to Brownian Motion, and accommodates for different drift parameters $\mu_i$ of $W_i(t)$. Karlin-McGregor's result can be applied to general strong Markov processes with continuous paths but no drifts. The transition probability for the process $W_i(t)$ is $\varphi_{s,\mu_i}(a-c)dc={\rm Pr}(W_i(s)\in dc \,|\, W_i(0)=a)$, where $\varphi_{s,\mu_i}(a-c)$ is the transition density.

\begin{corollary}
\be\label{key_BM_form}
&& {\rm Pr}\left\{W_0(t)<W_1(t)<\cdots<W_n(t), 0\le t\le s \, | \, W_i(0) =a_i,W_i(s)=c_i \,\,  (0\leq i \leq n)\right\} \nonumber \\
 &&\qquad\qquad\qquad = \exp\left(-\frac{s}{2}|\boldsymbol{\mu}|^2\!+\!\boldsymbol{\mu}\cdot(\textbf{c}-\textbf{a})  \right)\det\left[\varphi_{s}({a}_i-{c}_j) \right]_{i,j=0}^n\big/ \prod_{i=0}^{n}\varphi_s({a}_i\!-\!{c}_i\!+\!{\mu}_i).\;\;\;
\ee
\end{corollary}
\textbf{Proof}. Using the relation $\varphi_{s,\mu_i}(a-c)=\varphi_{s}(a-(c-\mu_i))$ and  dividing both sides of \eqref{Katori} by ${\rm Pr}(W_i(s)\in dc_i, i=0,1,\ldots,n\, | \, W_i(0)=a_i, i=0,1,\ldots,n)$, we
obtain the result.

\hfill $\Box$

\subsection{The main result}\label{main_result_1}
Let $\varphi(t)=\varphi_1(t)$
and $\Phi(t)=\int_{-\infty}^t \varphi(u)du$ be the density and the c.d.f. of the standard normal distribution.
Assume that $T=n$ is a positive integer. Define $(n\!+\!1)$-dimensional vectors
\begin{align}\label{bzy}
 \boldsymbol{\mu}  =
          \begin{bmatrix}
           0 \\
           b\\
           2b\\
           \vdots \\
           nb
         \end{bmatrix},
          \;\; \textbf{a}  =
          \begin{bmatrix}
           0 \\
           x_1\!+\!a\\
          x_2\!+\!2a\!+\!b\\
           \vdots \\
           x_n\!+\!na\!+\!\frac{(n-1)n}{2}b
         \end{bmatrix},
           \;\; \textbf{c}  =
          \begin{bmatrix}
           x_1 \\
           x_2\!+\!a\!+\!b\\
          x_3\!+\!2a\!+\!3b\\
           \vdots \\
           x_{n+1}\!+\!(a\!+\!b)n\!+\!\frac{(n-1)n}{2}b
         \end{bmatrix}
  \end{align}
and let ${\mu}_i$, ${a}_i$ and ${c}_i$ be $i$-th components of vectors $\boldsymbol{\mu}$, $\textbf{a}$ and $\textbf{c}$
respectively ($i=0,1, \ldots, n$). Note  that we start the indexation of vector components at 0.

\begin{theorem}
\label{th:1}
For any integer $n\geq 1 $ and $ x<a$,
\be\label{final_eqn}
\qquad F_{a,b}(n\, |\, x) \!=\!\frac{1}{\varphi(x)} \int_{-x-a-b}^{\infty}\int_{x_2-a-2b}^{\infty}\!\!\!\!\!\!\!\!\!&&\!\!\ldots \int_{x_n-a-nb}^{\infty}\exp(-|\boldsymbol{\mu}|^2/2 + \boldsymbol{\mu}\cdot(\textbf{c}-\textbf{a})) \nonumber\\
&&\qquad\quad   \times\det\left[\varphi({a}_i-{c}_j) \right]_{i,j=0}^n \, dx_{n+1}\,dx_{n}\!\ldots \,dx_{2}\,,\;\;\;\;
\ee
where $\boldsymbol{\mu}$, \textbf{a} and \textbf{c} are given in \eqref{bzy}.
\end{theorem}

If $b=0$ then \eqref{final_eqn} coincides with Shepp's formula (2.15) in \cite{Shepp71} expressed in variables $y_i=x_i+ia$ ($i=0,1,\ldots,n$).

%
%
In the case  $T=2$ we obtain
\bea
F_{a,b}(2\, |\, x)\!\!\!\!&=&\!\!\!\!\frac{e^{5b^2/2+bx}}{\varphi(x)}\int_{-x-a-b}^{\infty}\int_{x_2-a-2b}^{\infty} e^{b(2x_3-x_2)} \\&\times & \!\!\!\! \det\begin{bmatrix}
    \varphi(x)      &  \varphi(-x_2\!-\!a\!-\!b) &  \varphi(-x_3\!-\!2a\!-\!3b) \\
    \varphi(a)      & \varphi(-x\!-\!x_2\!-\!b)&   \varphi(-x\!-\!a\!-\!x_3\!-\!3b) \\
    \varphi(x_2\!+\!2a\!+\!b\!+\!x)      & \varphi(a)&  \varphi(x_2\!-\!x_3\!-\!2b) \\
\end{bmatrix} dx_3dx_2.
\eea

\subsection{An alternative representation of formula \eqref{final_eqn} and two particular cases}
It is easier to interpret Theorem~\ref{th:1} by expressing the integrals in terms of the values of $S(t)$ at  times $t=0,1, \ldots, n$.
Let
 $x_0= 0, x_1=-x.$ For $i=0,1, \ldots, n$ we set $s_i=x_i-x_{i+1}$ with $s_0=x$.
It follows from the proof of \eqref{final_eqn}, see Section~\ref{sec:proof}, that  $s_0, s_1, \ldots, s_n$ have the meaning of the values of the process $S(t)$ at  times $t=0,1, \ldots, n$; that is,
$S(i)=s_i$ ($i=0,1, \ldots, n$). The range of the variables $s_i$ in \eqref{final_eqn}  is $(-\infty,a+bi)$,  for $i=0,1, \ldots, n$ . The variables $x_1, \ldots, x_{n+1}$ are expressed via $s_0, \ldots, s_{n}$ by $x_k=-s_0-s_1-\ldots-s_{k-1}$
($k=1, \ldots, n+1$) with $x_0=0$. Changing the variables, we obtain the following equivalent expression for the probability $F_{a,b}(n\, |\, x)$:\\
\bea
F_{a,b}(n\, |\, x) \!=\!\frac{1}{\varphi(x)} \int_{-\infty}^{a+b}\int_{-\infty}^{a+2b}&&\!\!\!\!\!\!\!\!\!\!\! \cdots \int_{-\infty}^{a+bn}\exp(-|\boldsymbol{\mu}|^2/2 + \boldsymbol{\mu}\cdot(\textbf{c}-\textbf{a})) \nonumber\\
&&\qquad\!\! \times\det\left[\varphi({a}_i-{c}_j) \right]_{i,j=0}^n \, ds_{n}\!\ldots ds_{2}ds_{1}\, ,
\eea
where $\boldsymbol{\mu}$ is given by \eqref{bzy} but expressions for $\textbf{a}$ and $\textbf{c}$ change:
\begin{align*}
           \textbf{a}  =
          \begin{bmatrix}
           0 \\
           \!a\!-\!s_0\\
          2a\!+\!b \!-\! s_0\!-\!s_1\\
           \vdots \\
           na\!+\!\frac{(n-1)n}{2}b\!-\!s_0\!-\!s_1\!-\!\ldots\!-\!s_{n-1}
         \end{bmatrix}\!,\;\;
            \textbf{c}  =
          \begin{bmatrix}
           -s_0 \\
           a\!+\!b\!-\!s_0\!-\!s_1\\
          2a\!+\!3b\!-\!s_0\!-\!s_1\!-\!s_2\\
           \vdots \\
      (a\!+\!b)n\!+\!\frac{(n-1)n}{2}b\!-\!s_0\!-\!s_1\!-\!\ldots\!-\!s_{n}
         \end{bmatrix}.\;\;\;\;
  \end{align*}

In a particular  case of $n=1$ we obtain:
\be
\quad F_{a,b}(1\, |\, x)\!\!\!\!&=&\!\!\!\!\frac{1}{\varphi(x)}\int_{-\infty}^{a+b}\exp(-b^2/2\!+\!b(b-s_1))\det\begin{bmatrix}
    \varphi(x)      &  \varphi(x+s_1-a-b) \\
    \varphi(a)      & \varphi(s_1-b) \\
\end{bmatrix} ds_1 \nonumber \\
\!\!\!\!&=&\!\!\!\!\Phi(a+b)-\exp\left(-(a^2-x^2)/2-b(a-x) \right)\Phi(x+b), \label{T_1_expl}
\ee
which agrees with \eqref{case_T_1} in the Appendix A.
For the case of $n=2$ we obtain:

\bea
F_{a,b}(2\,|\,  x)\!\!\!\!&=&\!\!\!\!\frac{e^{5b^2\!/\!2 }}{\varphi(x)}\int_{-\infty}^{a\!+\!b}\! \int_{-\infty}^{a\!+\!2b} e^{ -b(s_1\!+\!2s_2)} \\
&\times &\!\!\!\! \det\!\!\begin{bmatrix}
    \varphi(x)    \!\!\!  & \!\!\! \varphi(x\!\!+\!\!s_1\!\!-\!\!a\!\!-\!\!b) &  \varphi(x\!\!+\!\!s_1\!\!+\!\!s_2\!\!-\!\!2a\!\!-\!\!3b) \\
    \varphi(a)   \!\!\!   & \!\!\varphi(s_1-b)&   \varphi(s_1+s_2-a-3b) \\
    \varphi(2a\!+\!b\!-\!s_1)  \!\!\!    & \varphi(a)&  \varphi(s_2-2b) \\
\end{bmatrix}\!\! ds_2 ds_1.
\eea

\subsection{Proof of Theorem~\ref{th:1}}
\label{sec:proof}

Using \eqref{BM} we rewrite $F_{a,b}(n\, |\, x)$  as
\bea
F_{a,b}(n\, |\, x)\!\!\!\! &=&\!\!\!\! {\rm Pr} \{W(t)\!-\!W(t\!+\!1)\!<a\!+\!bt \text{ for all   } t\!\in\![0,n]\,\, | \,\,W(0)\!-\!W(1)\!=\!x \} \\
\!\!\!\!&=&\!\!\!\! {\rm Pr} \{ W(t)\!-\!W(t\!+\!1)<a+bt,\,\, W(t\!+\!1)\!-\!W(t\!+\!2)<\!a+b(t\!+\!1),\ldots,\\
&&\quad \! W(t\!+\!n\!-\!1)\!-\!W(t\!+\!n)\!<\!a\!+\!b(t+\!n\!-1) \text{ for all   } t\!\in\![0,1] \,\, \big \lvert \,\,W(0)\!-\!W(1)\!=\!x \} \\
\!\!\!\!&=&\!\!\!\! {\rm Pr}\bigg\{W(t)\!<\!W(t\!+\!1)\!+\!a\!+\!bt<\!\cdots\!<W(t\!+\!n)\!+\!n(a\!+\!bt)\!+\!\frac{(n\!-\!1)n}{2}b\\
&& \quad \text{for all   } t\! \in\![0,1]\,\,
|\,\,W(0)\!-\!W(1)\!=\!x \bigg\}.
\eea
Let $\Omega$ be the event defined as follows
\bea
\Omega = \left\{ W(t)\!<\!W(t\!+\!1)\!+\!a\!+\!bt<\!\cdots\!<W(t\!+\!n)\!+\!n(a\!+\!bt)\!+\!\frac{(n\!-\!1)n}{2}b\text{ for all   } t\!\in\![0,1]  \right\}
\eea
and let $x_i=W(i)$, $i=0,1,\ldots,n\!+\!1$. Integrating out over the values $x_i$, by the law of total probability we obtain:
\be\label{Q_first}
\qquad F_{a,b}(n\, |\, x)\! = \!\int\!\! \!\!\!\! &\cdots &\!\!\!\!\!\!  \int {\rm Pr} \{ \Omega \,\, |\,\, W(0)\!=\!x_0,\ldots, W(n\!+\!1)\!=\!x_{n+1}, W(0)\!-\!W(1)\!=\!x \} \nonumber \\
& &\!\!\!\! \times \,{\rm Pr}\{ W(0) \!\in \!dx_0,\ldots, W(n\!+\!1) \!\in\! dx_{n+1} \, |\, W(0)\!-\!W(1)\!=\!x \}.
\ee
Note that $W(1)=x_1=-x$, since $W(0)-W(1)=x$ and $W(0)=0$. For $i=0,1,\ldots,n$, define the processes
\begin{align*}
W_i(t) = W(t+i)+i(a+bt)+ \frac{(i-1)i}{2}b, \,\,\,\, 0\le t \le 1.
\end{align*}
Then the event $\Omega$ above can be equivalently expressed as
\be\label{Omeg}
\Omega  = \{ W_0(t)<W_1(t)<\cdots < W_n(t) \text{ for all   } t\in[0,1]    \}
\ee
and under the conditioning introduced in \eqref{Q_first}, we have for $i=0,1,\ldots,n$:
\bea
W_i(0) \!\!\!&=&\!\!\! W(i)+ia + \frac{(i-1)i}{2}b = x_i + ia + \frac{(i-1)i}{2}b\, , \\
W_i(1) \!\!\!&=&\!\!\! W(i+1)+i(a+b)+\frac{(i-1)i}{2}b = x_{i+1}+i(a+b)+\frac{(i-1)i}{2}b \, .
\eea
Therefore \eqref{Q_first} can expressed as

\be\label{form}
&&F_{a,b}(n\, |\, x) \! =\!\! \int \!\!  \cdots\!\!\!\! \int \! {\rm Pr}\bigg\{ \Omega \, \big|\,  W_i(0)\!=\!x_i\!+\!ia\! +\!\frac{(i\!-\!1)i}{2}b, \; W_i(1)\!=\! x_{i+1}\!+\!i(a\!+\!b)\!+\!\frac{(i\!-\!1)i}{2}b\;(0 \leq i \leq n), \nonumber\\
& &    W_0(0)\!-\!  W_0(1)\!=\!x \! \bigg\} \; {\rm Pr} \{ W(0) \!\in\! dx_0, \ldots,\! W(n\!+\!1) \!\in\! dx_{n+1}
   \!  \,|\,W(0)\!-\!W(1)\!=\!x \}.\;\;\;\;\;\;\;\;\;\;\;\;\;\;\;\;\;\;
\ee

The region of integration for \eqref{form} is determined from the following chain of inequalities which ensure that the inequalities in \eqref{Omeg} hold at $t=0$ and $t=1$:
\bea
x_1<x_2\!+\!a\!+\!b<\ldots < x_n\!+\!(n-1)(a\!+\!b)\!+\!\frac{(n\!-\!2)(n\!-\!1)}{2}b< x_{n+1}\!+\!n(a\!+\!b)\!+\!\frac{(n\!-\!1)n}{2}b \,.
\eea
Hence, the upper limit of integration for all variables is infinity and the lower limit for the integral with respect to $x_{i+1}$, $i=1,\ldots,n$, is given by the formula:
\bea
x_i+(i-1)(a+b)+\frac{(i-2)(i-1)}{2}b-i(a+b) - \frac{(i-1)i}{2}b = x_i-a-b-(i-1)b \,.
\eea
Since the conditioned Brownian Motion processes $W_i(t)$ are independent, using \eqref{key_BM_form} we can express the first term in \eqref{form} as
\bea
{\rm Pr}\bigg\{ \!\Omega \,\big |\, W_i(0)\!\!\!\!&=&\!\!\!x_i\!+\!ia\! +\!\frac{(i\!-\!1)i}{2}b, W_i(1)\!=\! x_{i+1}\!+\!i(a\!+\!b)\!+\!\frac{(i\!-\!1)i}{2}b \text{  for  }i=0,1,\ldots n  \bigg \} \\
&=&\!\!\! \exp(-|\boldsymbol{\mu}|^2/2 + \boldsymbol{\mu}\cdot(\textbf{c}-\textbf{a}))\det\left[\varphi({a}_i-{c}_j) \right]_{i,j=0}^{n}\big /\prod_{i=0}^{n}\varphi({a}_i-{c}_i+{\mu}_i),
\eea
where $\boldsymbol{\mu}$, $\textbf{a}$ and $\textbf{c}$ are defined in \eqref{bzy}.
 The second probability in the right hand side of \eqref{form} is simply $\prod_{i=1}^{n}\varphi(x_i-x_{i+1})$.
By noticing
\bea
\prod_{i=0}^{n}\varphi({a}_i-{c}_i+{\mu}_i) = \prod_{i=0}^{n}\varphi(x_i-x_{i+1}-ib+{\mu}_{i})= \prod_{i=0}^{n}\varphi(x_i-x_{i+1})
\eea
and collating all terms, we obtain \eqref{final_eqn}.
\hfill $\Box$

\section{Linear barrier $a+bt$ with non-integral $T$}\label{sec:Non_int_T}

In this section, we shall derive an explicit formula for the first passage probability $F_{a,b}(T\, |\, x)$ defined in \eqref{problem} assuming $T>0$ is not an integer. Represent $T$ as $T=m+\theta$, where $m=\lfloor T \rfloor \ge 0$ is the integer part of $T$ and $0<\theta<1$. Set $n=m+1=\lceil T \rceil$.

\subsection{The main result}\label{main_result_2}

Let $\varphi_\theta(t)$ and $\varphi_{1-\theta}(t)$ be as defined  in \eqref{BM_transition}.
Define the $(n\!+\!1)$- and $n$-dimensional vectors as follows:  $\boldsymbol{\mu}_1= \boldsymbol{\mu}$ is as defined in \eqref{bzy},
\begin{align}\label{bzy_1}
           {\textbf{a}_1}  =
          \begin{bmatrix}
           0 \\
           u_1\!+\!a\\
          u_2\!+\!2a\!+\!b\\
           \vdots \\
           u_{n}\!+\!na\!+\!\frac{n(n-1)}{2}b
         \end{bmatrix},\;\;
            \textbf{c}_1  =
          \begin{bmatrix}
           v_0\\
           v_1\!+\!a\!+\!b\theta \\
          v_2\!+\!2(a\!+\!b\theta)\!+\!b\\
           \vdots \\
           v_{n}\!+\!n(a\!+\!b\theta)\!+\!\frac{n(n-1)}{2}b
         \end{bmatrix},
  \end{align}\;\;

\begin{align}
  \boldsymbol{\mu}_2  =
          \begin{bmatrix}
           0 \\
           b\\
           2b\\
           \vdots \\
           mb
         \end{bmatrix},
         \;  \textbf{a}_2  =
          \begin{bmatrix}
           v_0 \\
           v_1\!+\!a\!+\!b\theta\\
          v_2\!+\!2(a\!+\!b\theta)\!+\!b\\
           \vdots \\
           v_{m}\!+\!m(a\!+\!b\theta)\!+\!\frac{(m-1)m}{2}b
         \end{bmatrix},
           \; \textbf{c}_2  =
          \begin{bmatrix}
           u_1\\
           u_2\!+\!a\!+\!b \\
          u_3\!+\!2a\!+\!3b\\
           \vdots \\
           u_{m+1}\!+\!m(a\!+\!b)\!+\!\frac{(m-1)m}{2}b
         \end{bmatrix},\nonumber\\ \label{bzy_2}
  \end{align}
and let ${{a}_1}_i$ and ${{c}_1}_i$ be $i$-th components of vectors $\textbf{a}_1$ and $\textbf{c}_1$ respectively ($i\!=\!0,1, \ldots, n$). Similarly, let  ${{a}_2}_i$ and ${{c}_2}_i$ be $i$-th components of vectors  $\textbf{a}_2$ and $\textbf{c}_2$ respectively ($i\!=\! 0,1, \ldots, m$). Recall that we start the indexation of vector components at 0.

\begin{theorem}\label{th:2}
For  $ x<a$ and non-integral $T=m+\theta$ with $0<\theta<1$, we have
\bea
F_{a,b}(T\, |\, x) \!\!\!\!\!&=&\!\!\!\!\! \frac{1}{\varphi(x)} \int_{-x-a-b}^{\infty} \!\cdots\!  \int_{u_m-a-mb}^{\infty}\int_{-\infty}^{\infty}\int_{v_0-a-b\theta}^{\infty}\!\cdots  \int_{v_m-a-b\theta -mb}^{\infty}\\
&&\!\!\! \exp(-\theta|\boldsymbol{\mu}_1|^2/2 + \boldsymbol{\mu}_1\cdot(\textbf{c}_1-\textbf{a}_1))\exp(-(1\!-\!\theta)|\boldsymbol{\mu}_2|^2/2 + \boldsymbol{\mu}_2\cdot(\textbf{c}_2-\textbf{a}_2))\\
&&\!\!\! \det [\varphi_\theta({{a}_1}_i-{{c}_1}_j )]_{i,j=0}^n \det [\varphi_{1-\theta}({{a}_2}_i-{{c}_2}_j )]_{i,j=0}^m\,\,dv_{m+1}\ldots dv_1dv_0 du_{m+1}\ldots du_2\, .
\eea
\end{theorem}

Proof is given below in Section~\ref{sec:pr1}.

If $b=0$ then the above formula for $F_{a,b}(T\, |\, x)$ coincides with Shepp's formula (2.25) in \cite{Shepp71} expressed in variables $x_i=u_i+ia$ and $y_i=v_i+ia$ ($i=0,1,\ldots,n$).

\subsection{Two particular cases of Theorem~\ref{th:2}}

Taking $m=0$ and hence $T=\theta$ yields the following
\bea
F_{a,b}(\theta \, |\, x)\!\! \!&=&\!\! \!\frac{e^{\theta b^2/2}}{\varphi(x)}\int_{-\infty}^{\infty} \int_{v_0-a-b\theta }^{\infty} e^{ b(v_1+x)} \varphi_{1-\theta}(v_0\!+\!x)  \\
&\times &\!\! \det \begin{bmatrix}
     \varphi_\theta(\!-v_0)       &  \varphi_\theta(-v_1\!-\!a\!-\!b\theta)    \\
     \varphi_\theta(-x\!+\!a\!-\!v_0)       \!\!\!\!\!&  \varphi_\theta(-x\!-\!v_1\!-\!b\theta)  \\
\end{bmatrix}
dv_1dv_0
\eea
which agrees numerically with \eqref{case_T_1} in the Appendix for $T<1$.
Taking $T=1+\theta$ yields
\bea
F_{a,b}(1\!+\!\theta \, |\, x)\! \!\!\!\!&=&\!\!\!\!\!\frac{1}{\varphi(x)} \int_{-x-a-b}^{\infty}\int_{-\infty}^{\infty} \int_{v_0-a-b\theta }^{\infty} \int_{v_1-a-b\theta -b}^{\infty} \!\! \exp(-(1\!-\!\theta)b^2/2 \!+\!b(u_2\!-\!v_1\!+\!b(1\!-\!\theta))\\
&\times&\!\!\!\! \exp(-\theta(b^2\!+\!4b^2)/2 \!+\! b(v_1\!+\!b\theta\!+\!x)+\!2b(v_2\!-\!u_2\!+\!2b\theta))\\
&\times&\!\!\!\! \det \begin{bmatrix}
    \varphi_{1-\theta}(v_0\!+\!x)      &   \varphi_{1-\theta}(v_0\!-\!u_2\!-\!a\!-\!b)   \\
     \varphi_{1-\theta}(v_1\!+\!a\!+\!b\theta\!+\!x)      &  \varphi_{1-\theta}(v_1\!-\!u_2\!-\!b(1\!-\!\theta))   \\
\end{bmatrix}
\\
&\times&\!\! \!\!  \det \begin{bmatrix}
     \varphi_\theta(\!-v_0)       &  \varphi_\theta(-v_1\!-\!a\!-\!b\theta)   &   \varphi_\theta(-v_2\!-\!2(a\!+\!b\theta)\!-\!b)   \\
     \varphi_\theta(-x\!+\!a\!-\!v_0)       &  \varphi_\theta(-x\!-\!v_1\!-\!b\theta)  & \varphi_\theta(-x\!-\!v_2\!-\!a\!-\!b(1\!+\!2\theta))   \\
     \varphi_\theta(u_2\!+\!2a\!+\!b\!-\!v_0)       \!\!\! &  \varphi_\theta(u_2\!-\!v_1\!+\!a\!+\!b(1\!-\!\theta))  \!\!\!\!\!&  \varphi_\theta(u_2\!-\!v_2\!-\!2b\theta)   \\
\end{bmatrix}\\
&& \qquad \qquad \qquad \qquad \qquad \qquad\qquad\qquad\qquad\qquad\qquad\qquad\qquad dv_2dv_1dv_0 du_2 \, .
\eea

\subsection{Proof of Theorem~\ref{th:2}}
\label{sec:pr1}

We are interested in an expression for the first passage probability
\bea
F_{a,b}(T\, |\, x) = {\rm Pr}(S(t)<a+bt \text{ for all   } t\in[0,m+\theta]\,\, | \,\,S(0)=x).
\eea
Using \eqref{BM}, $F_{a,b}(T\, |\, x)$  can be equivalently expressed as follows
\bea
F_{a,b}(T\, |\, x)\!\!\!\!\! &=&\!\!\!\! {\rm Pr}\{ W(t)\!-\!W(t\!+\!1)\!<\!a\!+\!bt \text{ for all   } t\!\in\![0,m+\theta]\,\, | \,\,W(0)\!-\!W(1)\!=\!x\} \\
&=&\!\!\!\! {\rm Pr}\bigg\{ W(t)\!<\!W(t\!+\!1)\!+\!a\!+\!bt\!<\!\ldots \!<\!W(t\!+\!m\!+\!1)\!+\!(m\!+\!1)(a\!+\!bt)\!+\! \frac{(m\!+\!1)m}{2}b\\
&& \quad \text{for all   } t\!\in\![0,\theta] \text{  and  } W(\tau\!+\!\theta)\!<\!W(\tau\!+\!\theta\!+\!1)\!+\!a\!+\!b\theta\!+\!b\tau \!<\! \ldots \!<\! \\
&&  \quad W(\tau\!+\!\theta\!+\!m)+\!m(a\!+\!b\theta\!+\!b\tau)\!+\!\frac{(m\!-\!1)m}{2}b \text{ for all   } \tau \!\in\![0,1-\theta]  |\,W(0)\!-\!W(1)\!=\!x \bigg\} .
\eea
Let $\Omega$ be the event
\bea
\Omega \!&=&\! \bigg\{W(t)\!<\!W(t\!+\!1)\!+\!a\!+\!bt\!<\!\ldots\!<\!W(t\!+\!m\!+\!1)\!+\!(m\!+\!1)(a\!+\!bt)\!+\!\frac{(m\!+\!1)m}{2}b  \\
&& \quad\! \text{for all   } t\!\in\![0,\theta] \text{  and  } W(\tau\!+\!\theta)\!<\!W(\tau\!+\!\theta\!+\!1)\!+\!a\!+\!b\theta\!+\!b\tau\! < \!\ldots\! <\! W(\tau\!+\!\theta\!+\!m)+ \\
&&\quad \!   m(a\!+\!b\theta\!+\!b\tau)  +\frac{(m\!-\!1)m}{2}b \text{ for all   } \tau\!\in\![0,1-\theta] \bigg\}.
\eea
Then by integrating out over the values $u_i$ and $v_i$ of $W$ at times $i$ and $i+\theta$, $i=0,1,\ldots,m+\!1$, by the law of total probability we have
\be\label{Q_first_non_int}
F_{a,b}(T\, |\, x) =\! \int\!\!\!\!\! &\cdots& \!\!\!\!\! \int \!{\rm Pr}\{ \Omega \, \big| \,  W(0)\!=\!u_0,\ldots, \!W(m\!+\!1)\!=\!u_{m+1}, W(\theta)\!=\!v_0,\ldots,\! \nonumber \\
&& \qquad\!\!  W(m\!+\!1\!+\!\theta)\!=\!v_{m+1},W(0)\!-W(1)\!=\!x \} \times \nonumber
 \\
 && \, {\rm Pr}\{ W(0)\!\in\! du_0,\ldots, \!W(m\!+\!1) \!\in\! du_{m+1}, W(\theta)\!\in\! dv_0  ,\ldots,  \nonumber \\
 &&\qquad\!\!\! W(m\!+\!1\!+\!\theta) \!\in\! dv_{m+1} \, |\,W(0)\!-\!W(1)\!=\!x\}.
\ee
Note that $W(1)=x_1=-x$, since $W(0)-W(1)=x$ and $W(0)=0$, and define the processes
\begin{eqnarray*}
W_i(t)&=&W(t+i)+i(a+bt)+ \frac{(i-1)i}{2}b, \,\,\,\, 0\le t \le \theta, \,\,\,\,i=0,1,\ldots,m+1\,,\\
W_j^\prime(t)&=& W(\tau+\theta+j)+j(a+b\theta+b\tau)+\frac{(j-1)j}{2}b, \,\,\,\, 0\le \tau \le 1-\theta, \,\,\,\,j=0,1,\ldots,m\,.
\end{eqnarray*}
Then the event $\Omega$ can be equivalently expressed as $\Omega  = \Omega_1 \cap \Omega_2 $
with
\bea
\Omega_1\!\!\!&=&\!\!\!\{ W_0(t)<W_1(t)<\cdots < W_{m+1}(t) \text{ for all   } t\in[0,\theta]    \},\\
\Omega_2\!\!\!&=&\!\!\!\{ W_0^\prime(\tau)<W_1^\prime(\tau)<\cdots < W_m^\prime(\tau) \text{ for all   } \tau\in[0,1-\theta]    \}.
\eea
Under the  conditioning introduced in \eqref{Q_first_non_int} we have for $i\!=\!0,1,\ldots,m+1$ and $j\!=\!0,1,\ldots,m$:
\bea
W_i(0) \!\!\!&=&\!\!\! W(i)+ia + \frac{(i-1)i}{2}b = u_i + ia + \frac{(i-1)i}{2}b\, ,\\
W_i(\theta)\!\!\! &=&\!\!\! W(i+\theta)+i(a+b\theta)+\frac{(i-1)i}{2}b = v_{i}+i(a+b\theta)+\frac{(i-1)i}{2}b\, ,\\
W_j^\prime(0) \!\!\!&=&\!\!\! W(j+\theta)+j(a+b\theta)+\frac{(j-1)j}{2}b = v_j+j(a+b\theta)+\frac{(j-1)j}{2}b\, ,\\
W_j^\prime(1-\theta) \!\!\!&=&\!\!\! W(j+1)+j(a+b)+\frac{(j-1)j}{2}b = u_{j+1}+j(a+b)+\frac{(j-1)j}{2}b\, .
\eea
Now under the above conditioning the processes are independent and so the conditional probability of $\Omega$ in \eqref{Q_first_non_int} becomes a product of the conditional probabilities of $\Omega_1$ and $\Omega_2$. Therefore,  \eqref{Q_first_non_int} becomes
\be\label{split}
F_{a,b}(T\, |\, x) \!= \!\!\int\!\!\!\!\! &\cdots& \!\!\!\!\!\int \!{\rm Pr}\bigg\{ \!\Omega_1 \,\big|\,  W_i(0) \!=\!  u_i \!+\! ia \!+\! \frac{(i\!-\!1)i}{2}b, W_i(\theta) \!=\! v_{i}\!+\!i(a\!+\!b\theta)\!+\!\frac{(i\!-\!1)i}{2}b,\nonumber \\
&& \qquad     (0 \! \leq\! i\! \leq \! m\!+\!1) \bigg\} \times {\rm Pr}\bigg\{ \!\Omega_2\, \big| \, W_j^\prime(0)\!=\! v_j\!+\!j(a\!+\!b\theta)\!+\!\frac{(j\!-\!1)j}{2}b, \nonumber\\
&&\qquad  W_j^\prime(1\!-\!\theta)\!=\! u_{j+1}\!+\!j(a\!+\!b)\!+\!\frac{(j\!-\!1)j}{2}b \,\,  (0\! \leq\! j \leq\! m) \bigg\} \!\! \times \! {\rm Pr}\{ W(0)\! \in\! du_0 \nonumber \\
&& \qquad  ,\ldots,W(m\!+\!1) \!\in\! du_{m+1}, W(\theta) \!\in\! dv_0, \ldots,W(m\!+\!1\!+\!\theta) \! \in \!dv_{m+1} \nonumber\\
&& \qquad  |  W(0)\!-\!W(1)\!=\!x\}.
\ee
The region of integration for the variables $u_i$ in \eqref{split} is determined from the following chain of inequalities:
\bea
-x-a<u_2+2a+b<\ldots< u_m+ma+\frac{(m-1)m}{2}b<u_{m+1}+(m+1)a+\frac{(m+1)m}{2}b\,.
\eea
Whence, the upper limit of integration with respect to $u_{i+1}$ is infinity and the lower limit for the integral with respect to $u_{i+1}$, $i=1,\ldots,m$ is given by the formula $u_i-a-ib$. For the variables $v_j$ in \eqref{split}, we have the following chain of inequalities
\bea
v_0<v_1\!+\! a\!+\!b\theta \! <\ldots\!<v_{m}\!+\! m(a+b\theta)\!+\! \frac{(m-1)m}{2}b<v_{m+1}\!+\!(m+1)(a+b\theta)\!+\!\frac{(m+1)m}{2}b\,.
\eea
Once again, the upper limit of integration with respect to $v_{i+1}$ is infinity  and the lower limit for the integral with respect to $v_{i+1}$ ($i=0,\ldots,m$) is  $v_{i}-a-b\theta-ib$. For $v_0$, the upper and lower limits of integration are infinite.

Now using \eqref{key_BM_form} with $n=m+1$ we obtain
\bea
&&\!\!  {\rm Pr}\bigg\{ \! \Omega_1 \, \big |\, W_i(0)\! \!=\! \! u_i \!+\! ia \!+\! \frac{(i\!-\!1)i}{2}b, W_i(\theta) \!=\! v_{i}\!+\!i(a\!+\!b\theta)\!+\!\frac{(i\!-\!1)i}{2}b,\,\, (0\leq i \leq m\!+\!1) \bigg \} \\
&&\qquad = \exp(-\theta|\boldsymbol{\mu}_1|^2/2 + \boldsymbol{\mu}_1\cdot(\textbf{c}_1-\textbf{a}_1))\det [\varphi_\theta({{a}_1}_i-{{c}_1}_j )]_{i,j=0}^{m+1}\, \big /\prod_{i=0}^{m+1}\varphi_\theta({{a}_1}_i-{{c}_1}_i+\theta{{\mu}_1}_i),
\eea
where $\varphi_\theta(\cdot)$ is given in \eqref{BM_transition}, ${\textbf{a}_1}$  and  ${\textbf{c}_1}$ are given in \eqref{bzy_1}. Similarly, using \eqref{key_BM_form} with $n=m$ we have
\bea
 && \!\!\!\! {\rm Pr}\bigg\{ \!\Omega_2  \, \big| \,    W_j^\prime(0)\!=\! v_j\!+\!j(a\!+\!b\theta)\!+\!\frac{(j\!-\!1)j}{2}b,W_j^\prime(1\!-\!\theta)\!=\! u_{j+1}\!+\!j(a\!+\!b)\!+\!\frac{(j\!-\!1)j}{2}b, (0\! \leq\! j\! \leq\! m)\! \bigg \}  \\
&&\qquad = \exp(-(1\!-\! \theta)|\boldsymbol{\mu}_2|^2/2 + \boldsymbol{\mu}_2\cdot(\textbf{c}_2-\textbf{a}_2))\det [\varphi_{1-\theta}({{a}_2}_i-{{c}_2}_j )]_{i,j=0}^{m}\\
&&  \qquad \qquad\qquad\qquad\qquad\qquad\qquad\qquad \big /\prod_{i=0}^{m}\varphi_{1-\theta}({{a}_2}_i\! -\! {{c}_2}_i+(1\!-\! \theta){{\mu}_2}_i),
\eea
where
$\varphi_{1-\theta}(\cdot)$ is given in \eqref{BM_transition}, ${\textbf{a}_2}$  and  ${\textbf{c}_2}$ are given in \eqref{bzy_2}. The third probability in the right-hand side of \eqref{split} is simply
\bea
\frac{1}{\varphi(x)}\prod_{j=0}^{m}\prod_{i=0}^{m+1}\varphi_\theta(u_i-v_i)\varphi_{1-\theta}(v_j-u_{j+1})\, .
\eea
By noticing
\bea
\prod_{j=0}^{m}\prod_{i=0}^{m+1}\varphi_\theta({{a}_1}_i\!-\!{{c}_1}_i+\theta{{\mu}_1}_i)
\varphi_{1-\theta}({{a}_2}_j\!-\!{{c}_2}_j\!+\!(1\!-\!\theta){{\mu}_2}_j)=\prod_{j=0}^{m}\prod_{i=0}^{m+1}\varphi_\theta(u_i\!-\! v_i)\varphi_{1-\theta}(v_j\!-\!u_{j+1})
\eea
and collating all terms, we obtain the result.
\hfill $\Box$

\section{Piecewise linear barrier with one change of slope} \label{sec:Special_case_barrier}

\subsection{Formulation of the main result}

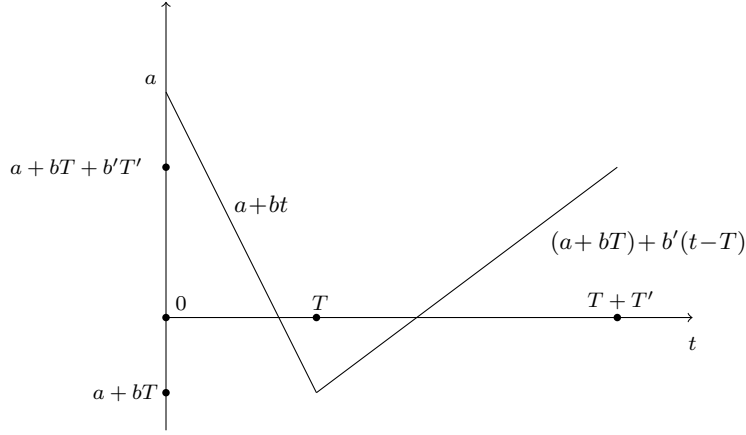
\begin{figure}[h]
\begin{center}
\begin{tikzpicture}

      \draw[->] (0,0) -- (7,0) node[right] {};
      \draw[->] (0,-1.5) -- (0,4.2) node[above] {};
      \draw[scale=1,domain=0:2,smooth,variable=\x]  plot ({\x},{3-2*\x});
      \draw[scale=1,domain=2:6,smooth,variable=\x]  plot ({\x},{-1+0.75*(\x-2)});
      \draw[scale=1,domain=4:4,smooth,variable=\x] plot ({\x},{3});

       \draw [-] (-0.2,3) node [above] {$a$};

       \draw [-] (-0.2,2) node [left] {$a+bT+ b^\prime T^\prime$};
 \node at (0,2) [circle,fill,inner sep=1pt]{};
       \draw [-] (2.05,0) node [above] {$T$};
       \draw [-] (6.05,0) node [above] {$T+T^\prime$};
       \draw [-] (0.2,0) node [above] {$0$};
       \draw [-] (0.8,1.5) node [right, font=\small] {$a\!+\!bt $};
       \draw [-] (5,1) node [right,font=\small] {$(a\!+bT)\!+b^\prime (t\!-\!T)$};
       \draw [-] (0,-1) node [left] {$a+bT$};
       \node at (2,0) [circle,fill,inner sep=1pt]{};
       \node at (6,0) [circle,fill,inner sep=1pt]{};
       \node at (0,0) [circle,fill,inner sep=1pt]{};
       \node at (0,-1) [circle,fill,inner sep=1pt]{};
        \draw [-] (7,-0.15) node [below] {$t$};

  \end{tikzpicture}
 \end{center}
 \vspace{-0.4cm}\caption{ Graphical depiction of a general boundary $B_{T,T^\prime}(t;a,b,b^\prime)$ with negative $b$ and positive $b'$. }\label{piecewise_barrier_bprime}
 \end{figure}

In this section, we derive an explicit formula for the first passage probability for $S(t)$ with a continuous piecewise linear barrier, where not more than  one change of slope is allowed.
For any non-negative $T,T'$ and real $a,b,b'$ we define the piecewise-linear barrier $B_{T,T^\prime}(t;a,b,b^\prime)$  by
\bea
B_{T,T^\prime}(t;a,b,b^\prime)= \begin{cases}
 a+bt & t \in [0,T], \\
 a+bT + b^\prime (t-T) & t \in [T, T+ T^\prime] \, ;
\end{cases}
\eea
for an illustration of this barrier, see Figure~\ref{piecewise_barrier_bprime}. We are interested in finding an expression for the first passage probability
\be\label{piecewise-problem3}
\qquad \qquad F_{a,b,b^\prime}(T,T^\prime \, |\, x) := {\rm Pr}(S(t)< B_{T,T^\prime}(t;a,b,b^\prime)\text{  for all   } t\in[0,T+T^\prime]\,\, | \,\,S(0)=x).
\ee
We only consider the case when both $T$ and $T'$ are integers. The case of general $T,T'$ can be treated similarly but the resulting expressions are much more complicated.

Define the $(T+T^\prime\!+\!1)$-dimensional vectors as follows:
\begin{align}\label{bzy3}
  \boldsymbol{\mu}_3  =
          \begin{bmatrix}
           0 \\
           b\\
           2b\\
           \vdots\\
           T b\\
           b^\prime  +T b\\
           2b^\prime  +T b\\
           \vdots\\
           T^\prime b^\prime  +Tb
         \end{bmatrix},
           \textbf{a}_3  =
          \begin{bmatrix}
           0 \\
           x_1+a\\
          x_2+2a+b\\
          \vdots \\
          x_{T}+Ta +\frac{(T-1)T}{2}b\\
          x_{T+1}+(T+1)a+bT +\frac{(T-1)T}{2}b \\
          x_{T+2}+(T+2)a+2bT +b^\prime+ \frac{(T-1)T}{2}b \\
          \vdots \\
         x_{T+T^\prime}+(T+T^\prime)a+bT T^\prime+\frac{(T^\prime-\!1)T^\prime}{2}b^\prime+ \frac{(T-1)T}{2}b
         \end{bmatrix},
  \end{align}

  \begin{align}\label{c_3}
            \textbf{c}_3  =
          \begin{bmatrix}
           x_1 \\
           x_2+a+b\\
          x_3+2a+3b\\
          \vdots\\
        x_{T}+(T-1)(a+b)+\frac{(T-2)(T-1)}{2}b\\
          x_{T+1}+T(a+b)+\frac{(T-1)T}{2}b\\
          x_{T+2}+a(T+1)+bT+ \frac{(T-1)T}{2}b +b^\prime +T b\\
          \vdots \\
          x_{T+T^\prime+1}+a(T+T^\prime)+bT T^\prime +\frac{(T^\prime-\!1)T^\prime}{2}b^\prime+ \frac{(T-1)T}{2}b+T^\prime b^\prime+Tb.
         \end{bmatrix},
  \end{align}
and let ${{a}_3}_i$ and ${{c}_3}_i$ be $i$-th components of vectors $\textbf{a}_3$ and $\textbf{c}_3$ respectively ($i=0,1, \ldots, T+T^\prime $).

\begin{theorem}\label{th:3}
For $ x<a$ and any positive integers $T$ and $T^\prime$, we have
\be\label{theorem3_form}
&& \!\! F_{a,b,b^\prime}(T,T'\, |\, x) \!=\!\frac{1}{\varphi(x)} \int_{-x-a-b}^{\infty}\int_{x_2-a-2b}^{\infty}\!\!\ldots \int_{x_T-a-T b}^{\infty} \int_{x_{T+1}-a-bT-b^\prime}^{\infty}\cdots \int_{x_{T+T^\prime}-a-bT - b^\prime T^\prime}^{\infty}\nonumber\\
&& \qquad \qquad \qquad\,\,\, \exp(-|\boldsymbol{\mu}_3|^2/2 + \boldsymbol{\mu}_3\cdot(\textbf{c}_3-\textbf{a}_3)) \det\left[\varphi({{a}_3}_i-{{c}_3}_j) \right]_{i,j=0}^{T+T^\prime} \, dx_{T+ T^\prime+1}\ldots dx_{2} \,.\;\;\;\;\;\;\;\;
\ee
\end{theorem}
Since the proof of Theorem~\ref{th:3} is similar to the proofs of Theorems~\ref{th:1} and \ref{th:2}, this proof is relegated to Appendix B, see Section~\ref{theorem3proof}.
Note that if $b=b'$ then \eqref{theorem3_form} reduces to  \eqref{final_eqn} with  $n=T+T'$.

\subsection{Two particular cases of Theorem~\ref{th:3}}
\label{sec:42}

Below we consider two particular cases of Theorem~\ref{th:3}; first, the barrier   is
$B_{1,1}(t;a,-b,b)$ with $b>0$; second, the barrier is  $B_{1,1}(t;a,0,-b')$ with $b'>0$. See Figures~\ref{two_cases} and \ref{three_cases} for a depiction of both barriers. As we shall demonstrate in Section~\ref{sec:Change_point}, these cases are important for problems of change-point detection.

\begin{figure}[ht]
\begin{minipage}[b]{0.45\linewidth}
\centering
\begin{tikzpicture}[scale=0.9] 
      \draw[->] (0,0) -- (5,0) node[right] {};
      \draw[->] (0,-1.5) -- (0,4.2) node[above] {};
      \draw[scale=1,domain=0:2,smooth,variable=\x]  plot ({\x},{3-2*\x});
      \draw[scale=1,domain=2:4,smooth,variable=\x]  plot ({\x},{-1+2*(\x-2)});
      \draw[scale=1,domain=4:4,smooth,variable=\x] plot ({\x},{3});
       \draw [-] (-0.2,3) node [above] {$a$};
       \draw [-] (2,0) node [above] {$1$};
       \draw [-] (4,0) node [above] {$ 2$};
       \draw [-] (0.2,0) node [above] {$0 $};
       \draw [-] (0.8,1.5) node [right, font=\small] {$a\!-\!bt $};
       \draw [-] (3.3,1.5) node [right,font=\small] {$a\!-\!2b\!+\!bt$};
       \draw [-] (0,-1) node [left] {$a\!-\!b$};
       \node at (2,0) [circle,fill,inner sep=1pt]{};
       \node at (4,0) [circle,fill,inner sep=1pt]{};
       \node at (0,0) [circle,fill,inner sep=1pt]{};
       \node at (0,-1) [circle,fill,inner sep=1pt]{};
        \draw [-] (5,-0.15) node [below] {$t$};
\end{tikzpicture}

\caption{ Barrier $B_{1,1}(t;a,-b,b)$ with $b>0$.}
\label{two_cases}
\end{minipage}
\hspace{0.5cm}
\begin{minipage}[b]{0.45\linewidth}
\centering
\begin{tikzpicture}[scale=0.9] 
      \draw[->] (0,0) -- (5,0) node[right] {};
      \draw[->] (0,-1.5) -- (0,4.2) node[above] {};
      \draw[scale=1,domain=0:2,smooth,variable=\x]  plot ({\x},{3});
      \draw[scale=1,domain=2:4,smooth,variable=\x]  plot ({\x},{3-2*(\x-2)});
       \draw [-] (-0.2,3) node [above] {$a$};
       \draw [-] (2,0) node [above] {$1$};
       \draw [-] (4,0) node [above] {$ 2$};
       \draw [-] (0.2,0) node [above] {$0 $};
       \draw [-] (3,1.5) node [right, font=\small] {$a\!+\!b'\!-\!b't $};
       \draw [-] (0,-1) node [left] {$a\!-\!b'$};
       \node at (2,0) [circle,fill,inner sep=1pt]{};
       \node at (4,0) [circle,fill,inner sep=1pt]{};
       \node at (0,0) [circle,fill,inner sep=1pt]{};
       \node at (0,-1) [circle,fill,inner sep=1pt]{};
        \draw [-] (5,-0.15) node [below] {$t$};
\end{tikzpicture}

\caption{Barrier $B_{1,1}(t;a,0,\!-b')$ with $b'\!>\!0$.}
\label{three_cases}
\end{minipage}\end{figure}

For the barrier $B_{1,1}(t;a,-b,b)$, an application of Theorem~\ref{th:3} yields
\bea
F_{a,-b,b}(1,1 \, |\, x) \!\!\!\!&=&\!\!\!\! \frac{e^{b^2/2}}{\varphi(x)} \int_{-x-a+b}^{\infty}\int_{x_2-a}^{\infty}e^{ - b(x_2+x)} \nonumber\\
 &\times&\!\!\!\! \det \begin{bmatrix}
    \varphi(x)      &  \varphi(\!-x_2\!-\!a\!+\!b) &  \varphi(\!-x_3\!-\!2a\!+\!b) \\
    \varphi(a)      & \varphi(-x\!-\!x_2\!+\!b)&  \varphi(\!-x\!-\!x_3\!-\!a\!+\!b) \\
    \varphi(x_2\!+\!2a\!-\!b\!+\!x)      & \varphi(a)   \!\!&  \varphi(x_2\!-\!x_3)
\end{bmatrix} dx_3dx_2.\nonumber
 \eea
\be
\qquad =\frac{e^{b^2/2-bx}}{\varphi(x)} \int_{-x-a+b}^{\infty} e^{- bx_2} \det \begin{bmatrix}
    \varphi(x)      &  \varphi(\!-x_2\!-\!a\!+\!b) &  \Phi(\!-x_2\!-\!a\!+\!b) \\
    \varphi(a)      & \varphi(\!-x\!-\!x_2\!+\!b)&  \Phi(\!-x\!-\!x_2\!+\!b) \\
    \varphi(x_2\!+\!2a\!-\!b\!+\!x)      & \varphi(a)&  \Phi(a) \\
\end{bmatrix} dx_2\,. \;\;\;\;\;
 \ee

For $B_{1,1}(t;a,0,-b')$, Theorem~\ref{th:3} provides:
\begin{eqnarray}\label{flat_down}
F_{a,0,-b'}(1,1 \, |\, x)\!\!\!\!&=&\!\!\!\! \frac{e^{{b'}^2}}{\varphi(x)} \int_{-x-a}^{\infty}\int_{x_2-a+b'}^{\infty}e^{-b'(x_3-x_2) } \nonumber\\
 &\times&\!\!\!\!  \det \begin{bmatrix}
    \varphi(x)      &  \varphi(-x_2-a) &  \varphi(-x_3-2a+b') \\
    \varphi(a)      & \varphi(-x-x_2)&  \varphi(-x-x_3-a+b') \\
    \varphi(x_2+2a+x)      & \varphi(a)&  \varphi(x_2-x_3+b')
\end{bmatrix} dx_3dx_2.
 \end{eqnarray}

\section{Piecewise linear barrier with two changes in slope}\label{Two_slope_changes}
\subsection{Formulation of the main result}

Theorem~\ref{th:3} can be generalized to the case when we have more than one change in slope. In the general case, the formulas for the first-passage probability  become very complicated; they are already rather heavy in the case of one change in slope.

In this section, we consider just one particular barrier with two changes in slope. For real $a,b,b',b''$, define the barrier $B(t;a,b,b^\prime,b'')$ as
\bea
B(t;a,b,b^\prime,b'')= \begin{cases}
 a+bt, &\;t \in [0,1], \\
 a+b + b^\prime (t-1), &\;t \in [1,2], \\
 a+b + b^\prime+ b''(t-2), &\;t \in [2,3]\, .
\end{cases}
\eea
As will be explained in Section~\ref{sec:Change_point}, the corresponding first-passage probability
\be\label{piecewise-problem3}
\quad F_{a,b,b^\prime,b''}(3 |\, x) := {\rm Pr}(S(t)< B(t;a,b,b^\prime,b'')\text{  for all   } t\in[0,3]\,\, | \,\,S(0)=x)
\ee
 is important  for some change-point detection problems.

Define the four-dimensional vectors as follows:
\begin{align}\label{bzy4}
  \boldsymbol{\mu}_4  =
          \begin{bmatrix}
           0 \\
           b\\
           b+b'\\
           b+b'+b''
         \end{bmatrix},\;\;
           \textbf{a}_4  =
          \begin{bmatrix}
           0 \\
           x_1+a\\
          x_2+2a+b\\
         x_3+3a+2b+b'
         \end{bmatrix},\;\;
            \textbf{c}_4  =
          \begin{bmatrix}
           x_1 \\
           x_2+a+b\\
          x_3+2a+2b+b'\\
          x_4+3a+3b+2b'+b''
         \end{bmatrix}
  \end{align}
and let ${{a}_4}_i$ and ${{c}_4}_i$ be $i$-th components of vectors $\textbf{a}_4$ and $\textbf{c}_4$ respectively ($i=0,1,2 ,3$).

\begin{proposition}\label{lw:2}
For $ x<a$ and any real $a,b,b'$ and $b''$
\be
\quad  \qquad F_{a,b,b^\prime,b''}(3\, |\, x)\!\! \!\!&=&\!\!\!\!\frac{1}{\varphi(x)} \int_{-x-a-b}^{\infty}\int_{x_2-a-b-b'}^{\infty}\int_{x_3-a-b-b'-b''}^{\infty} \nonumber\\
&&\!\!\!\! \exp(-|\boldsymbol{\mu}_4|^2/2 + \boldsymbol{\mu}_4\cdot(\textbf{c}_4-\textbf{a}_4)) \det\left[\varphi({{a}_4}_i-{{c}_4}_j) \right]_{i,j=0}^{3} dx_4 dx_3 dx_{2}. \label{lem:2}\,\;\;\;\;
\ee
\end{proposition}
For the proof of Proposition~\ref{lw:2}, see Section~\ref{Lemma3proof} in Appendix.

\subsection{A particular case of Proposition~\ref{lw:2}}
\label{sec:52}
\begin{figure}[h]
\begin{tikzpicture}[scale=1] 
      \draw[->] (0,0) -- (7,0) node[right] {};
      \draw[->] (0,-1.5) -- (0,4.2) node[above] {};
      \draw[scale=1,domain=0:2,smooth,variable=\x]  plot ({\x},{3});
      \draw[scale=1,domain=2:4,smooth,variable=\x]  plot ({\x},{3-2*(\x-2)});
       \draw[scale=1,domain=4:6,smooth,variable=\x]  plot ({\x},{-1+2*(\x-4)});
       \draw [-] (-0.2,3) node [above] {$h$};
       \draw [-] (2,0) node [above] {$1$};
       \draw [-] (4,0) node [above] {$ 2$};
        \draw [-] (6,0) node [above] {$ 3$};
       \draw [-] (0.2,0) node [above] {$0 $};
       \draw [-] (2.8,1.5) node [right, font=\small] {$h\!+\!\mu\!-\!\mu t $};
       \draw [-] (5.3,1.5) node [right, font=\small] {$h\!-\!3\mu\!+\!\mu t $};
       \draw [-] (0,-1) node [left] {$h\!-\!\mu$};
       \node at (2,0) [circle,fill,inner sep=1pt]{};
         \node at (6,0) [circle,fill,inner sep=1pt]{};
       \node at (4,0) [circle,fill,inner sep=1pt]{};
       \node at (0,0) [circle,fill,inner sep=1pt]{};
       \node at (0,-1) [circle,fill,inner sep=1pt]{};
        \draw [-] (7,-0.15) node [below] {$t$};
\end{tikzpicture}
\caption{Barrier $B(t;h,0,\!-\mu,\mu)$ with $\mu\!>\!0$.}
\label{flat_down_up_barrier}
\end{figure}
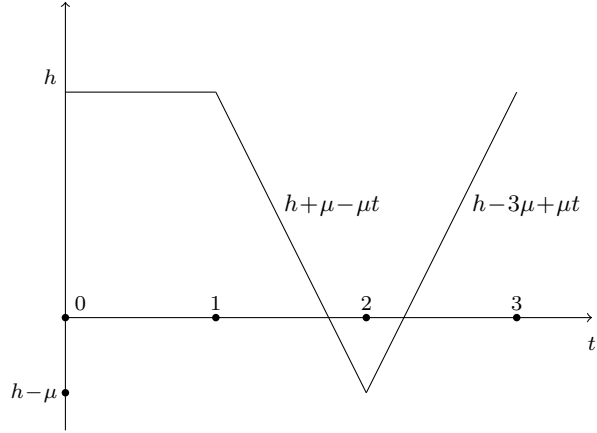

In this section, we consider a special barrier $B(t;h,0,-\mu,\mu)$ (depicted in Figure~\ref{flat_down_up_barrier}), which will be used in Section~\ref{sec:Change_point}. In the notation of Proposition~\ref{lw:2}, $a=h$, $b=0$, $b'=-\mu$, $b''=\mu$ and we obtain

\be \label{flat_down_up}
F_{h,0,-\mu,\mu}(3\, |\, x)=\frac{e^{{\mu}^2/2}}{\varphi(x)} \int_{-x-h}^{\infty}\int_{x_2-h+\mu}^{\infty} e^{ -\mu(x_3-x_2)} \times \;\;\;\;\;\;\;\;\;\;\;\;\;\;\;\;\;\;\;\;\;\;\;\;\;\;\;\;\;\;\;\;\;\;\;\;\;\;\;\;\;\;\; \;\;\; \;\;\; \;\;\; \;\;\nonumber \\
  \det \!\!  \begin{bmatrix}
    \varphi(x)      &  \varphi(-x_2\!-\!h) &  \varphi(-x_3\!-\!2h\!+\!\mu) &  \Phi(-x_3\!-\!2h\!+\!\mu) \\
    \varphi(h)      & \varphi(-x\!-\!x_2)&   \varphi(-x\!-\!x_3\!-\!h\!+\!\mu)\!\! & \Phi(-x\!-\!x_3\!-\!h\!+\!\mu) \\
    \varphi(x_2\!+\!2h\!+\!x)      & \varphi(h)&   \varphi(x_2\!-\!x_3\!+\!\mu)& \Phi(x_2\!-\!x_3\!+\!\mu) \\
        \varphi(x_3\!+\!3h\!-\!\mu\!+\!x)      \!\!& \varphi(x_3\!+\!2h\!-\!\mu\!-\!x_2) \!\!&   \varphi(h)& \Phi(h) \\
\end{bmatrix}\!\!dx_3dx_2.
\ee

\section{Application to change-point detection}\label{sec:Change_point}

In this section, we illustrate the natural appearance of first-passage probabilities for the Slepian process $S(t)$ for piece-wise barriers and in particular the barriers considered in Sections~\ref{sec:42} and \ref{sec:52}.

Suppose one can observe the stochastic process $X(t) \, (t\ge0)$ governed by the stochastic differential equation
\be
\label{eq:ch_p_W}
   dX(t) =\mu \mathbbm{1}_{\{\nu\leq t<\nu+l \}}dt + dW(t)\, ,
\ee
where $\nu>0$ is the unknown (non-random) change-point and  $\mu\neq0 $ is the drift magnitude during the `epidemic' period of duration $l$  with $0<l<\infty$; $\mu$ and $l$ may be known or unknown.
The classical  change-point detection problem of finding a change in drift of a Wiener process is the problem \eqref{eq:ch_p_W} with $l=\infty$; that is, when the change (if occurred) is permanent, see for example \citep{pollak1985diffusion,moustakides2004optimality,polunchenko2018asymptotic,polunchenko2010optimality}.

In \eqref{eq:ch_p_W}, under the null hypothesis $\mathbb{H}_0$, we assume $\nu = \infty$ meaning that the process $dX(t)$ has zero mean for all $t\ge 0$. On the other hand, under the alternative hypothesis $\mathbb{H}_1$,  $\nu<\infty$. In the definition of the test power, we will assume that $\nu$ is large. However, for the tests discussed below to be well-defined and approximations to be accurate, we only need $\nu \geq 1$ (under $\mathbb{H}_1$).

In this section, we only consider the case of known $l$, in which case we can assume $l=1$ (otherwise we change the time-scale by $t \to t/l$ and the barrier by $B \to B/\sqrt{l}$). The case when $l$ is unknown is more complicated and the  first-passage probabilities that have to be used are more involved; even so, these probabilities can be treated by the methodology similar to the one discussed below.

We define the test statistic used to monitor the epidemic alternative as
\bea
S_1(t) =\int_{t}^{t+1}dX(t)\,\,\,\, t\ge0\, .
\eea
The stopping rule for $S_{1}(t)$ is defined as follows
\be\label{tau_stopping}
\tau(h) = \inf \{t: S_{1}(t) \geq h   \},
\ee
where the threshold $h$ is chosen to satisfy the average run length (ARL) constraint $\mathbb{E}_0(\tau(h)  )=C$ for some (usually large) fixed $C$. Here $\mathbb{E}_0$ denote the expectation under the null  hypothesis.

For the process $S_1(t)- \mathbb{E}S_1(t)$, we  have
$$
S_1(t)- \mathbb{E}S_1(t)=W(t+1)-W(t)
$$
 which is stochastically equivalent to the Slepian process $S(t)$ of \eqref{BM}.

Under $\mathbb{H}_0$, $\mathbb{E}S_1(t) =0 $ for all $t\geq 0$ and under
$\mathbb{H}_1$ we have
$$
\mathbb{E}S_1(t)= \left\{
                    \begin{array}{cl}                      \mu(t-\nu) &\;\; {\rm for}\; \nu \leq t \leq \nu\!+\!1 \\
                      \mu -\mu(t-\nu-1)  & \;\;{\rm for}\; \nu+1 \leq t \leq \nu\!+\!2 \\
                      0 &\;\; {\rm otherwise}. \\
                    \end{array}
                  \right.
$$

The problem of construction of accurate approximations for $\mathbb{E}_0(\tau(h)  )$ relies on the  construction of accurate approximations  for the first-passage probabilities $ \int F_{h,0}(T\, |\, x) \varphi(x) dx$ for the Slepian process with constant barrier $h$ and large  $T$. This problem was addressed in \cite{noonan2018approximating}, where several accurate approximations were constructed. As a result, we can
derive an accurate approximation for  $\mathbb{E}_0(\tau(h)  )$, see  Section~\ref{ARL_section}. For example, to get $C =500$ we need $h\simeq 3.63$.
Since $l$ is known, for any $\mu>0$ the test with the stopping rule  \eqref{tau_stopping}
 is optimal in the sense  of the abstract Neyman-Pearson lemma, see Theorem 2, \cite[p 110]{grenander1981abstract}.

Here we are interested in the power of the test \eqref{tau_stopping} which can be defined as
 \be \label{Performance_measure}
\qquad {\cal P}({h,\mu}):= \lim_{\nu \to \infty}{\rm P}_1 \left\{S_{1}(t)\ge h \text{ for at least one } t\in[\nu,\nu+2 ]\,\, | \,\mathbb{H}_1,\; \tau(h)  > \nu \right\}\, ,
 \ee
where
 ${\rm P}_1$  denotes the probability measure  under the  alternative hypothesis.

 Define the piecewise linear barrier $Q_{\nu}(t;h,\mu) $ as follows
\bea
Q_{\nu}(t;h,\mu)  =h- \mu \max\{0,  1-|t-\nu - 1| \}.
\eea
The barrier $Q_{\nu}(t;h,\mu)$ is visually depicted in Fig~\ref{piecewise_barrier} below. The power of the test with the stopping rule \eqref{tau_stopping} is then
\bea
{\cal P}({h,\mu})=  \lim_{\nu \to \infty} {\rm P}\left\{S(t)\ge Q_{\nu}(t;h,\mu) \text{ for at least one } t\in[\nu,\nu+2]\,\, | \, \tau(h) > \nu \right\}.
\eea
Consider the barrier $B(t;h,0,-\mu,\mu)$ of Section~\ref{Two_slope_changes} with $t\in [0,3]$.
Define the conditional first-passage probability

\be
\gamma(x,h,\mu):= \nonumber
{\rm P}\{S(t)\ge B(t;h,0,-\mu,\mu) \text{ for at least one } t\in[1,3]\, \big| S(0)=x,\;   S(t)<h, \forall {t \in [0,1]} \} \nonumber
\\
=1- \frac{{\rm P}\left\{S(t)< B(t;h,0,-\mu,\mu) \text{ for all } t\in[0,3]\, \big| S(0)=x \right\}
} { {\rm P}\left\{S(t) <h \text{ for all } t\in[0,1] \big| S(0)=x \right\} } = 1- \frac{F_{h,0,-\mu,\mu}(3 |\, x)
} { F_{h,0}(1|x) }\, .\;\;\;\;\;\;
\label{eq:ratio}
\ee

The denominator in \eqref{eq:ratio}  is very simple to compute,  see \eqref{T_1_expl} with $b=0$ and $a=h$.
The numerator in \eqref{eq:ratio}  can be computed by \eqref{flat_down_up}.
Computation of  $\gamma(x,h,\mu)$ requires numerical evaluation of a two-dimensional integral, which is not  difficult.

We approximate the power ${\cal P}({h,\mu})$ by $\gamma(0,h,\mu)$.
In view of \eqref{BM} the process  $S(t)$ forgets the past after one unit of time hence quickly   reaches the stationary behaviour  under the condition $S(t)<h$ for all $t<\nu$. By approximating ${\cal P}({h,\mu})$ with $\gamma(0,h,\mu)$, we assume that one unit of time is almost enough for $S(t)$ to reach this stationary state. In Figure~\ref{fig:figure4}, we plot the ratio $\gamma(x,h,\mu)/\gamma(0,h,\mu)$ as a function of $x$ for $h=3$ and $\mu=3$. Since the ratio is very close to 1 for all considered $x$, this verifies that the probability $\gamma(x,h,\mu)$ changes very little as $x$ varies implying that the values of $S(t)$ at $T=\nu-1$ have almost no effect on the probability $\gamma(x,h,\mu)$. This allows us to claim that the accuracy $|{\cal P}({h,\mu}) -\gamma(0,h,\mu)|$ of the approximation
${\cal P}({h,\mu}) \simeq \gamma(0,h,\mu)$ is smaller than $10^{-5}$ for all $h \geq 3$. This claim agrees with discussions below in this section and  extensive simulations which we have performed. This claim also agrees with Table 2 in \cite{noonan2018approximating} (the row corresponding to $\lambda^{(4)}(h)$), from where we deduce that the accuracy of approximation ${\cal P}({h,\mu}) \simeq \gamma(0,h,\mu)$ is smaller than $10^{-6}$ for all $h \geq 3$ and $\mu=0$; it is also intuitively clear that the accuracy of the  approximation  ${\cal P}({h,\mu}) \simeq \gamma(0,h,\mu)$ improves as $\mu $ grows.

\begin{figure}[h]
\begin{minipage}[b]{0.4\linewidth}
\centering
\raisebox{0.5cm}{
\begin{tikzpicture}[scale=1]

      \draw[->] (-1.5,0) -- (6,0) node[right] {};
      \draw[->] (0,-1.5) -- (0,3.7) node[above] {};
      \draw[scale=1,domain=0:2,smooth,variable=\x]  plot ({\x},{3-2*\x});
      \draw[scale=1,domain=2:4,smooth,variable=\x]  plot ({\x},{-1+2*(\x-2)});
      \draw[scale=1,domain=4:4,smooth,variable=\x] plot ({\x},{3});
      \draw[scale=1,domain=4:6,smooth,variable=\x] plot ({\x},{3});
      \draw[scale=1,domain=-1.5:0,smooth,variable=\x] plot ({\x},{3});
       \draw [-] (-0.2,3) node [above] {$h$};
       \draw [-] (2,0) node [above] {$\nu\! +\! 1$};
       \draw [-] (4,0) node [above] {$\nu \!+ \!2$};
       \draw [-] (0.2,0) node [above] {$\nu $};
       \draw [-] (0.8,1.5) node [right, font=\small] {$\!h\!\!-\!\!\mu( t\!\!-\!\!\nu) $};
       \draw [-] (3.3,1.5) node [right,font=\small] {$\!h\!\!-\!\!\mu \!+\!\mu (t\!\!-\!\!\nu\!\!-\!\!1)$};
       \draw [-] (0,-1) node [left] {$h\!-\!\mu $};
       \node at (2,0) [circle,fill,inner sep=1pt]{};
       \node at (4,0) [circle,fill,inner sep=1pt]{};
       \node at (0,0) [circle,fill,inner sep=1pt]{};
       \node at (0,-1) [circle,fill,inner sep=1pt]{};
           \node at (0,3) [circle,fill,inner sep=1pt]{};
        \draw [-] (6,-0.15) node [below] {$t$};
  \end{tikzpicture}}

\caption{Graphical depiction of the {boundary} $Q_{\nu}(t;h,\mu)$.}
\label{piecewise_barrier}
\end{minipage}
\hspace{2cm}
\begin{minipage}[b]{0.45\linewidth}
\centering
    \includegraphics[scale=0.18]{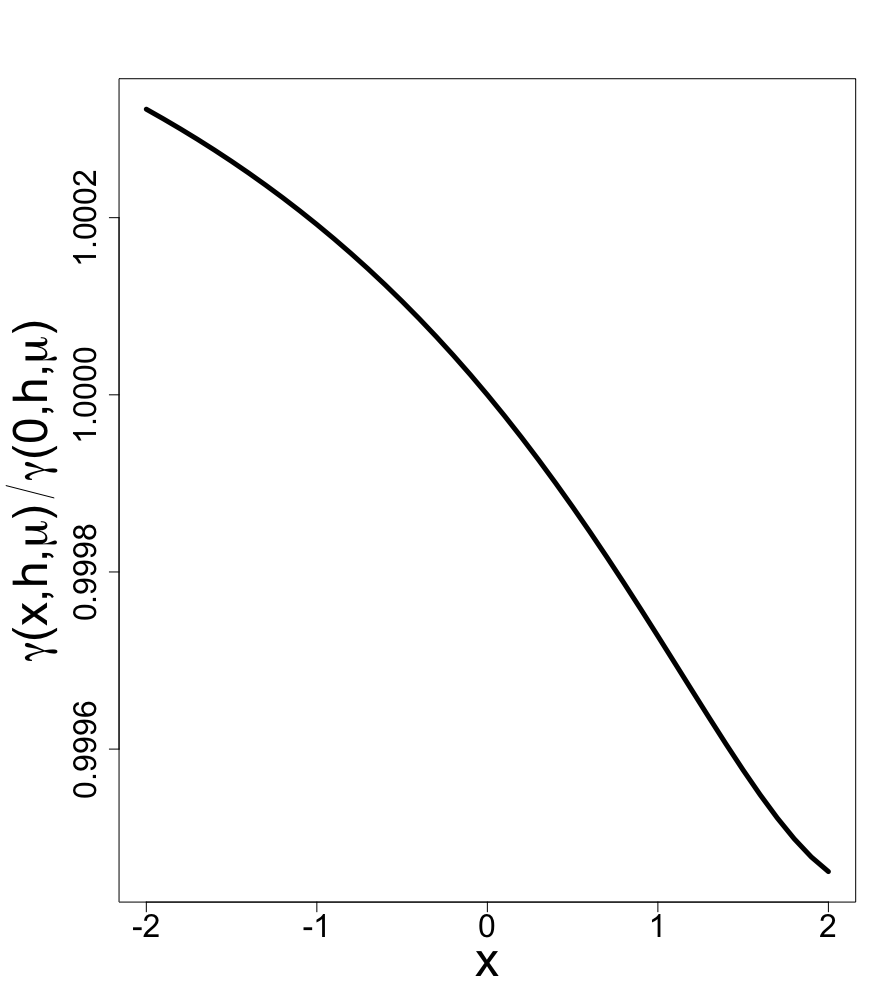}
\caption{Ratio $\gamma(x,h,\mu)/\gamma(0,h,\mu)$ for $h=3$ and $\mu=3$.}
\label{fig:figure4}
\end{minipage}
\end{figure}

 In Table~\ref{power_table}, we provide values of  $\gamma(0,h,\mu)$ for different $\mu$, where the values of $h$ have been chosen to satisfy $\mathbb{E}_0(\tau(h)  )=C$ for $C=100,500,1000$; see  \eqref{ARL_form_app} regarding computation of the ARL $\mathbb{E}_0(\tau(h)  )$.

 As seen from Figures~\ref{two_cases} and \ref{flat_down_up_barrier}, the barrier $B_{1,1}(t;h,-\mu,\mu)$  is the main component of the barrier $B(t;h,0,-\mu,\mu)$. Instead of using the approximation ${\cal P}({h,\mu}) \simeq \gamma(0,h,\mu)$ it is therefore tempting to use a simpler approximation
 ${\cal P}({h,\mu}) \simeq \gamma_1(0,h,\mu)$, where
\bea
\gamma_1(x,h,\mu):=
{\rm P}\{S(t)\ge B(t;h,-\mu,\mu) \text{ for at least one } t\in[0,2]\, \big| S(0)=x \}
= 1- F_{h,-\mu,\mu}(1,1 |\, x)
\eea
To compute  values of $\gamma_1(0,h,\mu)$ we only need to  evaluate a one-dimensional integral. Table~\ref{power_table_2}
 we show some values of  $\gamma_1(0,h,\mu)$ for different $\mu$. Comparing the entries of Tables \ref{power_table} and \ref{power_table_2} we can observe that the quality of
 ${\cal P}({h,\mu}) \simeq \gamma_1(0,h,\mu)$ is not too bad, especially for large $\mu$.

 Approximation ${\cal P}({h,\mu}) \simeq \gamma_1(0,h,\mu)$ can be improved if we average values of $\gamma_1(x,h,\mu)$ over an appropriate distribution for $x$. According to Section 2.4.2 in \cite{noonan2018approximating}, one of possible appropriate distributions for $x$ has density
  \bea
 p(x)=\frac{\Phi(h)\varphi(x) - \Phi(x)\varphi(h)}{ \Phi^2(h)-\varphi(h) \big[h\Phi(h)
+ \varphi(h) \big]}\, , \;\; x \leq h\, .
 \eea
 Define $\gamma_2(h,\mu)= \int_{-\infty}^h \gamma_1(x,h,\mu)p(x) dx $, which is a two-dimensional integral. As seen from comparison of Tables~\ref{power_table} and \ref{power_table_3}, the accuracy of the approximation ${\cal P}({h,\mu}) \simeq \gamma_2(h,\mu)$ is almost the same as the accuracy of the main approximation
 ${\cal P}({h,\mu}) \simeq \gamma(0,h,\mu)$. Computational cost of computing $\gamma_2(h,\mu)$ is similar to the cost for $\gamma(h,\mu) $.

 To assess the impact the final line-segment in the barrier  $B(t;h,0,-\mu,\mu)$ on power (the line-segment with gradient $\mu$ in Fig~\ref{piecewise_barrier}, $t \in [\nu,\nu+1]$), in Table~\ref{power_table_4} we document the values of
 $\gamma_3(0,h,\mu)$ for different $\mu$. Here
\bea
\gamma_3(x,h,\mu)\!\!\!&:=&\!\!\! \nonumber
{\rm P}\{S(t)\ge B(t;h,0,-\mu) \text{ for at least one } t\in[1,2]\, \big| S(0)=x,\; \,   S(t)<h, \forall {t \in [0,1]} \} \nonumber
\\
&=&\!\!\!1- \frac{{\rm P}\left\{S(t)< B(t;h,0,-\mu) \text{ for all } t\in[0,2]\, \big| S(0)=x \right\}
} { {\rm P}\left\{S(t) <h \text{ for all } t\in[0,1] \big| S(0)=x \right\} } = 1- \frac{F_{h,0,-\mu}(1,1 |\, x)
} { F_{h,0}(1|x) }
\eea
and  $F_{h,0,-\mu}(1,1  |\, 0)$ can be  computed using  \eqref{flat_down} with $b'=\mu$. By comparing Tables~\ref{power_table} and \ref{power_table_4}, one can see the expected diminishing impact which the final line-segment in $B(t;h,0,-\mu,\mu)$ has on power, as $\mu$ increases. However, for small $\mu$ the contribution of this part of the barrier to power is significant suggesting it is not be sensible to approximate the power of our test with $\gamma_3(0,h,\mu)$.

\begin{table}[h]
\begin{center}
\begin{tabular}{ |c|c|c|c| }
  \hline
  \multicolumn{2}{|c|}{ $  h=3.11$, $C\simeq 100$ } \\\hline
  $\mu$ & $\gamma(0,3.11,\mu)$\\ \hline
  2 &0.3052  \\
   2.25 &0.3876  \\
  2.5 &  0.4765\\
  2.75 & 0.5676  \\
  3 & 0.6559\\
     3.25 &0.7371 \\
  3.5 &0.8075  \\
  3.75 & 0.8653 \\
  4 &  0.9101\\
        4.25&0.9429  \\
  4.5 & 0.9655  \\
  4.75 & 0.9802 \\
  5 &0.9892 \\
  \hline
\end{tabular}
\begin{tabular}{ |c|c|c|c| }
  \hline
  \multicolumn{2}{|c|}{$  h=3.63$, $C\simeq 500$} \\\hline
  $\mu$ & $\gamma(0,3.63,\mu)$\\ \hline
  2 &0.1384  \\
   2.25 & 0.1946 \\
  2.5 & 0.2638 \\
  2.75 &0.3445  \\
  3 &0.4338 \\
       3.25 & 0.5272 \\
  3.5 &  0.6197\\
  3.75 & 0.7061 \\
  4 &0.7824 \\
          4.25&0.8461  \\
  4.5 & 0.8961 \\
  4.75 &  0.9332 \\
  5 & 0.9592\\
  \hline
\end{tabular}
\begin{tabular}{ |c|c|c|c| }
  \hline
  \multicolumn{2}{|c|}{ $  h=3.83$, $C\simeq 1000$} \\\hline
  $\mu$ & $\gamma(0,3.83,\mu)$\\ \hline
  2 &  0.0956  \\
   2.25 &0.1402  \\
  2.5 & 0.1979 \\
  2.75 &0.2687  \\
  3 &0.3510 \\
       3.25 &0.4416  \\
  3.5 &0.5358  \\
  3.75 &0.6285  \\
  4 &0.7146 \\
          4.25&0.7900  \\
  4.5 & 0.8525  \\
  4.75 & 0.9011 \\
  5 & 0.9370\\
  \hline
\end{tabular}
\end{center}
\caption{$\gamma(0,h,\mu)$ as a function of $\mu$ for three choices of ARL.}
\label{power_table}
\end{table}

\begin{table}[h]
\begin{center}
\begin{tabular}{ |c|c|c|c| }
  \hline
  \multicolumn{2}{|c|}{ $  h=3.11$, $C\simeq 100$ } \\\hline
  $\mu$ & ${\gamma_1}(0,3.11,\mu)$\\ \hline
  2 & 0.2918   \\
  2.5 & 0.4645  \\
  3 &  0.6471   \\
  3.5 & 0.8021  \\
  4 &  0.9075    \\
  4.5 &0.9644 \\
  5 &     0.9889  \\
  \hline
\end{tabular}
\begin{tabular}{ |c|c|c|c| }
  \hline
  \multicolumn{2}{|c|}{$  h=3.63$, $C\simeq 500$} \\\hline
  $\mu$ & ${\gamma_1}(0,3.63,\mu)$\\ \hline
  2 &   0.1310    \\
  2.5 &    0.2553   \\
  3 &  0.4256  \\
  3.5 &  0.6132   \\
  4 &  0.7783  \\
  4.5 &  0.8940    \\
  5 &  0.9583   \\
  \hline
\end{tabular}
\begin{tabular}{ |c|c|c|c| }
  \hline
  \multicolumn{2}{|c|}{ $  h=3.83$, $C\simeq 1000$} \\\hline
  $\mu$ & ${\gamma_1}(0,3.83,\mu)$\\ \hline
  2 & 0.0903    \\
  2.5 &  0.1911   \\
  3 &    0.3438 \\
  3.5 &  0.5295 \\
  4 &  0.7101  \\
  4.5 &   0.8499   \\
  5 &   0.9358  \\
  \hline
\end{tabular}
\end{center}
\caption{Values of $\gamma_1(0,h,\mu)$ for some  $\mu$ and $h$.}
\label{power_table_2}
\end{table}

\begin{table}[!h]
\begin{center}
\begin{tabular}{ |c|c|c|c| }
  \hline
  \multicolumn{2}{|c|}{ $  h=3.11$, $C\simeq 100$ } \\\hline
  $\mu$ & ${\gamma_2}(3.11,\mu)$\\ \hline
  2 & 0.3047    \\
  2.5 &0.4760   \\
  3 &   0.6555   \\
  3.5 &  0.8073  \\
  4 &   0.9100   \\
  4.5 & 0.9654 \\
  5 &   0.9892  \\
  \hline
\end{tabular}
\begin{tabular}{ |c|c|c|c| }
  \hline
  \multicolumn{2}{|c|}{$  h=3.63$, $C\simeq 500$} \\\hline
  $\mu$ & ${\gamma_2}(3.63,\mu)$\\ \hline
  2 &  0.1383   \\
  2.5 &    0.2637    \\
  3 &    0.4337  \\
  3.5 & 0.6196    \\
  4 &  0.7824 \\
  4.5 &   0.8961    \\
  5 &   0.9592  \\
  \hline
\end{tabular}
\begin{tabular}{ |c|c|c|c| }
  \hline
  \multicolumn{2}{|c|}{ $  h=3.83$, $C\simeq 1000$} \\\hline
  $\mu$ & ${\gamma_2}(3.83,\mu)$\\ \hline
  2 &  0.0956     \\
  2.5 &0.1978      \\
  3 & 0.3509    \\
  3.5 & 0.5358    \\
  4 & 0.7146   \\
  4.5 &   0.8524      \\
  5 &   0.9370  \\
  \hline
\end{tabular}
\end{center}
\caption{ Values of $\gamma_2(h,\mu)$ for some  $\mu$ and $h$.}
\label{power_table_3}
\end{table}

\begin{table}[!h]
\begin{center}
\begin{tabular}{ |c|c|c|c| }
  \hline
  \multicolumn{2}{|c|}{ $  h=3.11$, $C\simeq 100$ } \\\hline
  $\mu$ & ${\gamma_3}(0,3.11,\mu)$\\ \hline
  2 & 0.2389  \\
  2.5 & 0.4017 \\
  3 &   0.5873\\
  3.5 & 0.7567 \\
  4 &  0.8801 \\
  4.5 & 0.9514  \\
  5 &  0.9840 \\
  \hline
\end{tabular}
\begin{tabular}{ |c|c|c|c| }
  \hline
  \multicolumn{2}{|c|}{$  h=3.63$, $C\simeq 500$} \\\hline
  $\mu$ & ${\gamma_3}(0,3.63,\mu)$\\ \hline
  2 &0.1039   \\
  2.5 & 0.2131   \\
  3 & 0.3731 \\
  3.5 & 0.5611   \\
  4 & 0.7373  \\
  4.5 &  0.8685  \\
  5 & 0.9458  \\
  \hline
\end{tabular}
\begin{tabular}{ |c|c|c|c| }
  \hline
  \multicolumn{2}{|c|}{ $  h=3.83$, $C\simeq 1000$} \\\hline
  $\mu$ & ${\gamma_3}(0,3.83,\mu)$\\ \hline
  2 & 0.0708    \\
  2.5 & 0.1575  \\
  3 &  0.2974 \\
  3.5 & 0.4785   \\
  4 &  0.6657 \\
  4.5 & 0.8194     \\
  5 &  0.9192 \\
  \hline
\end{tabular}
\end{center}
\caption{Values of $\gamma_3(0,h,\mu)$ for some  $\mu$ and $h$. }
\label{power_table_4}
\end{table}

\section{Appendix A}
\label{sec:appA}
\subsection{First-passage probability $F_{a,b}(T\, |\, x)$ for $T\leq 1$}
\label{sec:T1}
For $T\le1$, the first passage probability $F_{a,b}(T\, |\, x)$ has been well studied. An explicit formula was first derived in 1988 in \cite[p.81]{ZhK1988} (published in Russian) and more than 20 years later it was independently derived in \cite{bischoff2016boundary} and \cite{deng2017boundary}. The authors of \cite{ZhK1988} and \cite{deng2017boundary} also considered the case of piecewise-linear barriers.

In \cite{ZhK1988}, the first passage probability $F_{a,b}(T\, |\, x)$ for $T\le 1$ was obtained by using the fact  $S(t)$ is a conditionally Markov process on the interval $[0,1]$. It was shown in \cite{Mehr} that after conditioning on $S(0)=x$, $S(t)$ can be expressed in terms of Brownian Motion as follows
\bea
S(t)=(2-t)W(g(t))+x(1-t), \,\,\,\,\,\,\,\,  0\le t\le 1.
\eea
 with $g(t)=t/(2-t)$. From this it follows that for $T\le 1$
\begin{align*}
F_{a,b}(T\,|\,x)={\rm Pr}\bigg ( W(g(t)) < \frac{a+bt-x(1-t)}{2-t} \text{ for all  } t\in[0,T]  \bigg ).
\end{align*}
Noting that   $t={2g(t)}/({1+g(t)})$ and using the well known barrier crossing formula for the Brownian motion  (see e.g.  \cite{Sieg_paper})
\begin{align}
F_{a,b}(T\,|\,x)=& \, {\rm Pr}\bigg ( W(g(t)) < \bigg( \frac{(a-x)(1+g(t))}{2} \bigg) + (x+b)g(t) \text{ for all } t\in[0,T]  \bigg) \nonumber\\
=& \, {\rm Pr}\bigg(W(t^\prime) < \bigg(\frac{a-x}{2}\bigg) + t^\prime\bigg(\frac{a+x}{2} + b\bigg) \text{ for all } t^\prime \in \bigg[0,\frac{T}{2-T} \bigg] \bigg ) \nonumber \\
 =&\, \Phi\bigg(\frac{b_1Z+a_1}{\sqrt{Z}}\bigg) - e^{-2a_1b_1}\Phi\bigg(\frac{b_1Z-a_1}{\sqrt{Z}}\bigg)\, \label{case_T_1}\, ,
\end{align}
where $Z ={T}/({2-T})$, $b_1= ({a+x})/{2} + b $ and $a_1 = ({a-x})/{2}$. This methodology, like many others, fails for $T>1$.



\subsection{An approximation for $\mathbb{E}_0(\tau(h)  )$} \label{ARL_section}
Consider the unconditional probability  (taken with respect to the standard normal distribution):
\bea
F_{h,0}(T) := \int_{-\infty}^{h} F_{h,0}(T\, |\, x)\varphi(x) dx \,.
\eea
Under  $\mathbb{H}_0$, the distribution of $\tau(h) $
  has the form:
  \bea
  (1-\Phi(h)) \delta_0 (d s) + q_h(s) d s \, ,s \geq 0,
  \eea
   where $\delta_0 (d s)$ is the delta-measure concentrated at 0 and
    \begin{equation*}
q_h(s)=-\frac{d}{ds}F_{h,0}(s ) ,\;\;\;0<s<\infty\, 
\end{equation*}
is the first-passage density. 
This yields
\be\label{ARL_form}
\mathbb{E}_0(\tau(h) ) = \int_{0}^{\infty}sq_h(s) ds.
\ee
There is no easy computationally convenient  formula for $q_h(t)$ as expressions for $F_{h,0}(s )$ are very complex.
 For deriving approximations for $\mathbb{E}_0(\tau(h) )$ we apply approximations for $F_{h,0}(s )$, discussed in \cite{noonan2018approximating}. One of the simplest (yet  very accurate) approximation takes the following form:
\be\label{F_2_approx}
F_{h,0}(T )\simeq F_{h,0}(2 ) \cdot \lambda(h)^{T-2}, \;\;  \mbox{for all $T>0$},
\ee
with $\lambda(h) = F_{h,0}(2)/  F_{h,0}(1)$. Values of $F_{h,0}(2 )$ must be numerically computed; approximations and simpler forms of $F_{h,0}(2 )$ have been presented in  \cite{noonan2018approximating} should one require an explicit formula. Using \eqref{F_2_approx}, we approximate the  density  $q_h(s)$ by
\bea
 q_h(s) \simeq  - F_{h,0}(2)\log [\lambda({h})] \cdot \lambda({h})^{s -2 }, \,\,\, 0< s <\infty.
\eea
Evaluation of the integral in \eqref{ARL_form} yields 
\be\label{ARL_form_app}
\mathbb{E}_0(\tau(h)  ) \cong -\frac{F_{h,0}(2 )}{ \lambda({h})^2\log [ \lambda({h})]}\,.
\ee
Numerical study shows that the approximation \eqref{ARL_form_app} is very accurate for all $h\geq 3$.

\section{Appendix B}
\label{sec:appB}

\subsection{Proof of \eqref{theorem3_form}}\label{theorem3proof}
The proof of  \eqref{theorem3_form} follows similar steps to the proof of \eqref{final_eqn}. The event~$\Omega$ becomes
\bea
\Omega &=& \bigg\{W(t)<W(t+1)+a+bt<\ldots < W(t+T)+ T(a+bt)+\frac{(T-1)T}{2}b \\
&&\,\,\,\,\,\, <W(t+T+1)+ a(T+1)+bT +\frac{(T-1)T}{2}b+(b^\prime+Tb)t \,\,\,\,\, < \ldots<\\
&&\quad  W(t+T +T^\prime)+a(T+T^\prime) +bT T^\prime + \frac{(T^\prime-1)T^\prime}{2}b^\prime + \frac{(T-1)T}{2}b +(T^\prime b^\prime +T b)t \\
&& \quad \text{for all   } t\in[0,1]   \bigg \}.
\eea
As in the proof of \eqref{final_eqn}, let $x_i=W(i)$, $i=0,\ldots,T+T^\prime +1 $, where $W(1)=x_1=-x$. Then
\be\label{F_first}
F_{a,b,b^\prime}(T,T'\, |\, x) \!=\! \int\!\! \! \!\!\!&\cdots&\!\!\!\! \!\!\int\! {\rm Pr}\{ \Omega\, | \, W(0)\!=\!x_0,\ldots, W(T+T^\prime +1)\!=\!x_{T+T^\prime +1}, W(0)\!-\!W(1)\!=\!x \} \nonumber\\
&&\!\! \!\! \times {\rm Pr}\{ W(0)\!\in\! dx_0,\ldots, W(T+T^\prime +1) \!\in\! dx_{T+T^\prime +1} \,\, |\,\,W(0)\!-\!W(1)\!=\!x\}.\;\;\;\;\;
\ee
Define the following processes which take different forms depending on the value of $i$:
\bea
W_i(t)\! &=&\! W(t+i)+i(a+bt)+\frac{(i-1)i}{2}b \,,\;\;\mbox{for $0\le i\leq T$}\,;
\\
W_i(t)\! &=&\!  W(t+i)\!+\!ai\!+\!bT (i-T) \!+\! \frac{(i\!-\!T\!-\!1)(i\!-\!T)}{2}b^\prime \!+\! \frac{(T\!-\!1)T}{2}b\!+\!\{(i\!-\!T)b^\prime\!+\!T b\}t ,
\eea
for  $T \!+\!1\leq i \leq T+ T^\prime\!$,
with $0\le t \le1$ for all processes. The event $\Omega$ can now be expressed as
\be\label{Omeg_f}
\qquad \Omega = \{W_0(t)<W_1(t)<\ldots<W_{T}(t)<\ldots <W_{T+T^\prime}(t) \text{ for all   } t\in[0,1]   \}.
\ee
Under the conditioning introduced in \eqref{F_first}, depending on the size of $i$ we have: for $ 0\le i \le T$
\bea
W_i(0) = x_i + ia+\frac{(i-1)i}{2}b\, ,  \;\;\;
W_i(1) = x_{i+1} + i(a+b)+\frac{(i-1)i}{2}b\, ;
\eea
and for $T +1 \leq i \leq T + T^\prime$
\bea
W_i(0) &=& x_i+ai +bT(i-T)\!+\! \frac{(i\!-\!T\!-\!1)(i\!-\!T)}{2}b^\prime \!+\! \frac{(T\!-\!1)T}{2}b\, , \\
W_i(1) &=& x_{i+1}+ai +bT(i-T)\!+\! \frac{(i\!-\!T\!-\!1)(i\!-\!T)}{2}b^\prime \!+\! \frac{(T\!-\!1)T}{2}b+\!(i\!-\!T)b^\prime\!+\!T b\, .
\eea
Whence \eqref{F_first} can be expressed as
\be\label{form2}
F_{a,b,b^\prime}(T,T'\, |\, x)\! =\! \int\!\! \!\!\! &\cdots&\!\!\! \!\!\int {\rm Pr} \bigg\{\Omega \, \bigg | \, W_i(0) = x_i + ia+\frac{(i-1)i}{2}b\,,W_i(1) = x_{i+1}  \nonumber\\
&&\qquad +\, i(a+b)+\frac{(i-1)i}{2}b\, \,\, (0\leq i \leq T), W_i(0) = x_i\!+\!ai\! +\!bT(i\!-\!T)  \nonumber \\
&& \qquad + \,  \frac{(i\!-\!T\!-\!1)(i\!-\!T)}{2}b^\prime \!+\! \frac{(T\!-\!1)T}{2}b\,, W_i(1) = x_{i+1}\!+\!ai \!+\! bT(i\!-\!T)\! \nonumber\\
&&\qquad + \,  \frac{(i\!-\!T\!-\!1)(i\!-\!T)}{2}b^\prime \!+\! \frac{(T\!-\!1)T}{2}b+\!(i\!-\!T)b^\prime\!+\!T b\,\,\, \nonumber \\
 && \qquad  (T\leq i \leq T\!+\!T^\prime), \,  W_0(0)\!-\!W_0(1)\!=\!x \bigg\}\times  \nonumber\\
&&{\rm Pr}\{ W(0)\!\in\! dx_0,\ldots, W(T\!+\!T^\prime \!+\!1) \!\in\! dx_{T+T^\prime +1} \, |\,W(0)\!-\!W(1)\!=\!x \}.\;\;\;\;
\ee
The region of integration in \eqref{form2} is determined from the following inequalities which, like in the proof of \eqref{theorem3_form}, ensure that the inequalities in \eqref{Omeg_f} hold at $t=0$ and $t=1$:
\bea
x_1\!\!\! &<\ldots <&\!\!\!  x_{T+1}\!+\!T (a+b) \!+\! \frac{(T\!-\!1)T}{2}b < x_{T+2}\!+\!a(T\!+\!1)\!+\!bT\!+\! \frac{(T\!-\!1)T}{2}b \!+\! b^\prime \!+\!T b<\ldots<\\
&& \!\!\! x_{T+T^\prime+1}+a(T+T^\prime)+bT T^\prime +\frac{(T^\prime\!-\!1)T^\prime}{2}b^\prime+ \frac{(T\!-\!1)T}{2}b +\!T^\prime b^\prime\!+\! T b.
\eea
From this, the upper limit of integration is infinity for all $x_i$. For $0\leq i \leq T +1$, the lower limit for  $x_{i}$ is  $x_{i-1}-a-(i-1)b$. For $T+2 \leq i \leq T +T^\prime+1$, the lower limit for  $x_i$ is  $x_{i-1}-a-bT - b^\prime (i-T-1)$.

Application of  \eqref{key_BM_form} with $n=T+T'$  provides
\bea
{\rm Pr} \bigg\{\!\!\!\!\!\! &\Omega&\!\!\!\!\!\! \, \bigg | \, W_i(0) = x_i + ia+\frac{(i-1)i}{2}b\,,W_i(1) = x_{i+1} + i(a+b)+\frac{(i-1)i}{2}b \, \, (0\leq i \leq T) \nonumber\\
&&\!\! W_i(0) = x_i+ai +bT(i-T)\!+\! \frac{(i\!-\!T\!-\!1)(i\!-\!T)}{2}b^\prime \!+\! \frac{(T\!-\!1)T}{2}b\, , \nonumber \\
&& \!\! W_i(1) = x_{i+1}+ai +bT(i-T)\!+\! \frac{(i\!-\!T\!-\!1)(i\!-\!T)}{2}b^\prime \!+\! \frac{(T\!-\!1)T}{2}b+\!(i\!-\!T)b^\prime\!+\!T b\, \nonumber\\
&&\!\!  (T\leq i \leq T+T^\prime), W_0(0)\!-\!W_0(1)\!=\!x \bigg\} \nonumber\\
&=& \exp(-| {\boldsymbol{\mu}}_3|^2/2 +  {\boldsymbol{\mu}}_3\cdot(\textbf{c}_3-\textbf{a}_3))\det [\varphi({{a}_3}_i,{{c}_3}_j )]_{i,j=0}^{T^\prime+T}/\prod_{i=0}^{T+T^\prime}\varphi({{a}_3}_i-{{c}_3}_i+{{\mu}_3}_i),
\eea
where
 $\boldsymbol{\mu}_3$ and $\textbf{a}_3$ are given in \eqref{bzy3} and $\textbf{c}_3$ is given in \eqref{c_3}. The second probability in the right hand side of \eqref{form2} is $\prod_{i=1}^{T+T^\prime}\varphi(x_i-x_{i+1})$.
We complete the proof by collating all terms and noting 
\bea
\prod_{i=0}^{T+T^\prime}\varphi({{a}_3}_i-{{c}_3}_i+{{\mu}_3}_i) = \prod_{i=0}^{T +T^\prime}\varphi(x_i-x_{i+1})\, .
\eea
\hfill $\Box$


\subsection{Proof of \eqref{lem:2}}\label{Lemma3proof}

Like the proof of \eqref{theorem3_form}, the proof of  \eqref{lem:2} is similar to the proof of \eqref{final_eqn}. We modify the event $\Omega$ as follows:
\bea
\Omega &=& \bigg\{W(t)<W(t+1)+a+bt<W(t+2)+2a+b+bt+b't <\\
&&\quad W(t+3)+3a+2b+b' +(b+b'+b'')t\text{ for all   } t\in[0,1]   \bigg \}.
\eea
By the law of total probability,
\be\label{F_first_2}
F_{a,b,b^\prime,b''}(3\, |\, x) \!=\! \int\!\! \!\!\!&\cdots&\!\!\! \!\!\int {\rm Pr}\{ \Omega\, | \, W(0)\!=\!x_0,\ldots, W(4)\!=\!x_{4}, W(0)\!-\!W(1)\!=\!x \}  \nonumber\\
&&\!\!\! \times \, {\rm Pr}\{ W(0)\!\in\! dx_0,\ldots, W(4) \!\in\! dx_{4} \,\, |\,\,W(0)\!-\!W(1)\!=\!x \}.
\ee
Define individually the following processes:
\bea
W_0(t) &=& W(t)\\
W_1(t) &=& a+bt+W(t+1)\\
W_2(t)&=&2a+b+(b+b')t+W(t+2)\\
W_3(t) &=& 3a+2b+b'+(b+b'+b'')t+W(t+3)
\eea
with $0\le t \le1$ for all processes. The event $\Omega$ can be re-written as
\be\label{Omeg_f2}
\Omega = \{W_0(t)<W_1(t)<W_2(t)<W_{3}(t) \text{ for all   } t\in[0,1]   \}.
\ee

The conditioning introduced in \eqref{F_first_2} results in:

\noindent\begin{minipage}{.5\linewidth}
\bea
W_0(0) &=& 0\\
W_1(0) &=& a+x_1\\
W_2(0)&=&2a+b+x_2\\
W_3(0) &=& 3a+2b+b'+x_3
\eea
\end{minipage}%
\begin{minipage}{.5\linewidth}
\bea
W_0(1) &=& x_1\\
W_1(1) &=& a+b+x_2\\
W_2(1)&=&2a+2b+b'+x_3\\
W_3(1) &=& 3a+3b+2b'+b''+x_4.
\eea
\end{minipage} \\ \\

From this, we can express \eqref{F_first_2} as:
\be\label{form5}
F_{a,b,b^\prime,b''}(3\, |\, x)\! \!=\!\! \int\!\! \!\!\! &\cdots& \!\!\! \!\!\int {\rm Pr} \bigg\{\Omega \, \bigg | \, W_0(0) = 0,\ldots,W_3(0) \!=\! 3a\!+\!2b\!+\!b'\!+\!x_3 , W_0(1) = x_1,\ldots, \nonumber \\
&& \qquad W_3(1) = 3a+3b+2b'+b''+x_4,  W_0(0)\!-\!W_0(1)\!=\!x \bigg\}  \nonumber\\
&&\,\,\,\, \times  {\rm Pr}\{ W(0)\!\in\! dx_0,\ldots, W(4) \!\in\! dx_{4} \,\, |\,\,W(0)\!-\!W(1)\!=\!x \}.
\ee

The region of integration for \eqref{form5} is determined from the following inequalities (see proof of \eqref{final_eqn} for similar discussion):
\bea
x_1<x_{2}+a+b  < x_{3}+2a+2b+b'<x_{4}+3a+3b+2b'+b''.
\eea
Thus, the upper limit of integration is infinity for all $x_i$. For integration with respect to $x_4$, the lower limit is $x_3-a-b-b'-b''$. For  integration with respect $x_3$,  the lower limit is  $x_2-a-b-b'$. Finally, for $x_2$, the lower limit  is $x_1-a-b=-x-a-b$. Now using \eqref{key_BM_form} with $n=3$ we obtain
\bea
{\rm Pr} \bigg\{\Omega \!\!\!\!\! & | & \!\!\!\!\!  W_0(0) = 0,\ldots,W_3(0) = 3a+2b+b'+x_3  \nonumber \\
&& \!\!\!\!\! W_0(1) = x_1,\ldots,W_3(1) = 3a+3b+2b'+b''+x_4,  W_0(0)\!-\!W_0(1)\!=\!x \bigg\} \nonumber\\
&=& \exp(-| {\boldsymbol{\mu}}_4|^2/2 +  {\boldsymbol{\mu}}_4\cdot(\textbf{c}_4-\textbf{a}_4))\det[\varphi({{a}_4}_i,{{c}_4}_j )]_{i,j=0}^3/\prod_{i=0}^{3}\varphi({{a}_4}_i-{{c}_4}_i+{{\mu}_4}_i),
\eea
 $\boldsymbol{\mu}_4$, $\textbf{a}_4$ and $\textbf{c}_4$ are given in \eqref{bzy4}. The second probability in the right hand side of \eqref{form5} is $\prod_{i=1}^{3}\varphi(x_i-x_{i+1})$.
Using the fact
\bea
\prod_{i=0}^{3}\varphi({{a}_4}_i-{{c}_4}_i+{{\mu}_4}_i) = \prod_{i=0}^{3}\varphi(x_i-x_{i+1}),
\eea
and collecting all results we complete the proof.
\hfill $\Box$

\bibliographystyle{unsrt}

\bibliography{changepoint}

\begin{thebibliography}{10}

\bibitem{slepian1961first}
D.~Slepian.
\newblock First passage time for a particular {G}aussian process.
\newblock {\em The Annals of Mathematical Statistics}, 32(2):610--612, 1961.

\bibitem{Mehr}
C.B. Mehr and J.A. McFadden.
\newblock Certain properties of {G}aussian processes and their first-passage
  times.
\newblock {\em Journal of the Royal Statistical Society. Series B
  (Methodological)}, 27(3):505--522, 1965.

\bibitem{Shepp71}
L.~Shepp.
\newblock First passage time for a particular {G}aussian process.
\newblock {\em The Annals of Mathematical Statistics}, 42(3):946--951, 1971.

\bibitem{karlin1959coincidence}
S.~Karlin and J.~McGregor.
\newblock Coincidence probabilities.
\newblock {\em Pacific Journal of Mathematics}, 9(4):1141--1164, 1959.

\bibitem{katori2012reciprocal}
M.~Katori.
\newblock Reciprocal time relation of noncolliding {B}rownian motion with
  drift.
\newblock {\em Journal of Statistical Physics}, 148(1):38--52, 2012.

\bibitem{pollak1985diffusion}
M.~Pollak and D.~Siegmund.
\newblock A diffusion process and its applications to detecting a change in the
  drift of {B}rownian motion.
\newblock {\em Biometrika}, 72(2):267--280, 1985.

\bibitem{moustakides2004optimality}
G.~Moustakides.
\newblock {O}ptimality of the {CUSUM} procedure in continuous time.
\newblock {\em The Annals of Statistics}, 32(1):302--315, 2004.

\bibitem{polunchenko2018asymptotic}
A.~Polunchenko.
\newblock Asymptotic near-minimaxity of the randomized
  {S}hiryaev--{R}oberts--{P}ollak change-point detection procedure in
  continuous time.
\newblock {\em Theory of {P}robability \& Its {A}pplications}, 62(4):617--631,
  2018.

\bibitem{polunchenko2010optimality}
A.~Polunchenko and A.~Tartakovsky.
\newblock On optimality of the {S}hiryaev--{R}oberts procedure for detecting a
  change in distribution.
\newblock {\em The Annals of Statistics}, 38(6):3445--3457, 2010.

\bibitem{noonan2018approximating}
J.~Noonan and A.~Zhigljavsky.
\newblock Approximating {S}hepp's constants for the {S}lepian process.
\newblock {\em arXiv preprint arXiv:1812.11101}, 2018.

\bibitem{grenander1981abstract}
U.~Grenander.
\newblock {\em Abstract inference}.
\newblock John Wiley \& Sons, 1981.

\bibitem{ZhK1988}
A.~Zhigljavsky and A.~Kraskovsky.
\newblock {\em Detection of abrupt changes of random processes in radiotechnics
  problems}.
\newblock St. Petersburg University Press, 1988.
\newblock (in Russian).

\bibitem{bischoff2016boundary}
W.~Bischoff and A.~Gegg.
\newblock Boundary crossing probabilities for (q, d)-{S}lepian-processes.
\newblock {\em Statistics \& Probability Letters}, 118:139--144, 2016.

\bibitem{deng2017boundary}
P.~Deng.
\newblock Boundary non-crossing probabilities for {S}lepian process.
\newblock {\em Statistics \& Probability Letters}, 122:28--35, 2017.

\bibitem{Sieg_paper}
D.~Siegmund.
\newblock Boundary crossing probabilities and statistical applications.
\newblock {\em The Annals of Statistics}, 14(2):361--404, 1986.

\end{thebibliography}


\end{document}